\documentclass[reqno,11pt]{amsart}

\newtheorem{theorem}{Theorem}[section]
\newtheorem{lemma}[theorem]{Lemma}
\newtheorem{corollary}[theorem]{Corollary} 

\theoremstyle{definition}

\theoremstyle{remark}
\newtheorem{remark}[theorem]{Remark}

\makeatletter
\def\dashint{\operatorname%
{\,\,\text{\bf--}\kern-.98em\DOTSI\intop\ilimits@\!\!}}
\makeatother

\newcommand{\nlimsup}{\operatornamewithlimits{\overline{lim}}}

\newcommand{\tr}{\text{\rm tr}\,}
\newcommand{\dist}{\text{\rm dist}\,}

\newcommand\bR{\mathbb{R}}

\newcommand\cF{\mathcal{F}}

\newcommand\sfu{{\sf u}}

\newcommand{\loc}{{\rm loc}}

\newcommand{\WO}{\overset{\scriptscriptstyle0}%
{W}\,\!}

 \newcommand{\mysection}[1]{\section{#1}
 \setcounter{equation}{0}}

\begin{document}

\title[Weighted Aleksandrov estimates]
{Weighted Aleksandrov estimates: PDE and
stochastic versions}
\author[N. V. Krylov]{N.V. Krylov}
\address[N. V. Krylov]{127 Vincent Hall, University of Minnesota,
 Minneapolis, MN, 55455}
\email{nkrylov@umn.edu}

\begin{abstract}We prove several pointwise  estimates
for solutions  of linear elliptic (parabolic) equations
with measurable coefficients in smooth domains (cylinders)
through the weighted $L_{d}$ ($L_{d+1}$)-norm of the free term.
The weights allow the free term to blow up near the 
(latteral) boundary. We also present weighted estimates
for occupation times of diffusion processes.
 
\end{abstract}

\maketitle

In the recent paper \cite{DK_18} the authors
prove weighted and mixed-norm $L_{p}$ estimates
for fully nonlinear elliptic and parabolic equations
with relaxed convexity assumption and almost
VMO assumption on the dependence on the space-time
variables  of the functions defining the equations.
The full norm including the second order spacial derivatives
is estimated, however,  sometimes they are estimated
   through the weighted norm of the free term
plus the weighted norm of  the unknown
function itself. Here we show how the unknown
function can be estimated  through the weighted norm of 
the free term by considering {\em linear\/}
elliptic and parabolic equations with {\em measurable\/} coefficients.
One knows that zeroth-order estimates for
fully nonlinear elliptic and parabolic equations even not explicitly
involving $x$ reduce  to the estimates for
linear elliptic and parabolic equations with measurable coefficients.
Our estimates are given for $C^{1,1}$ domains and cylinders
with weights that are powers of the distance
to the boundary of the domain or
 to the lateral boundary of the cylinder.
In \cite{Na_01} and the references therein one  can find similar
estimates in case the boundaries of domains
have wedges with weights related to the wedges.

It is worth noting that for the case of linear and quasilinear 
elliptic equations
with {\em regular\/} coefficients a rather detailed
information about weighted estimates
can be found in \cite{AN_02}
and references therein.

It is also worth noting that for the case of linear  
 parabolic equations
with   coefficients {\em independent\/} of $x$
a rather detailed
information about weighted estimates
can be found in \cite{KN_14}
and the references therein. It is also worth
noting \cite{Am_16}, where an abstract treatment
of weighted estimates is presented
from the point of view of semigroups
and special Riemannian manifolds.
Then the coefficient of operators are necessarily
smooth apart from some special singularities.

Our method in the elliptic case
 is an extension of the original
Aleksandrov methods based on Monge-Amp\`ere
equations and is presented in Section 
\ref{section 8.11.1}. In Section
\ref{section 8.11.2} we apply the results
of Section 
\ref{section 8.11.1} to derive estimates
for equations of main type in the unit ball.
Section \ref{section 8.15.1} contains our main
analytic result about estimates of solutions
of elliptic equations. 
In Section \ref{section 9.1.1}
we derive stochastic Aleksandrov estimates
for functions which can blow up near the 
boundary. These provide better estimates
than known before for the time spent by
diffusion processes near the boundary
of a domain before reaching it.

Our method in the parabolic case is an extension of the
one introduced in \cite{Kr_76},
 is based on considering the
parabolic
Monge-Amp\`ere
equations introduced in \cite{Kr_76}, and is presented in Section 
\ref{section 8.11.10} where we also
derive estimates for equations of main type in   round cylinders. In Section
\ref{section 8.11.20} we apply the results
of Section 
\ref{section 8.11.10} to derive estimates
for general parabolic 
equations   in   round cylinders.
Section \ref{section 8.15.10} contains our main
analytic result about estimates of solutions
of parabolic equations. Finally,
in Section \ref{section 9.1.10}
we derive stochastic weighted Aleksandrov estimates
for functions which can blow up near the 
boundary. These provide better estimates
than known before for the time spent by
diffusion processes near the 
lateral boundary
of a cylinder before reaching its boundary.

In the elliptic part of the article
we work in a $d$-dimensional
Euclidean space $\bR^{d}$ of points $ x=(x^{1},...,
x^{d}) $, $d\geq2$.
We use the notation
$$
D_{i}=\frac{\partial}{\partial x^{i}},
\quad Du=(D_{1}u,...,D_{d}u),\quad
D_{ij}=D_{i}D_{j},\quad D^{2}u=(D_{ij}u)\big|_{i,j=1}^{d},
$$
$$
a_{\pm}=(1/2)(|a|\pm a),\quad a_{\pm}^{p}=(a_{\pm})^{p},\quad
B_{r}(x)=\{y\in\bR^{d}:|y-x|<r\},
$$
$$
B_{r}=B_{r}(0),\quad B=B_{1}.
$$

By $|\Gamma|$ we denote the volume of $\Gamma\subset\bR^{d}$.
If $\Omega$ is a domain in $\bR^{d}$ with regular boundary
by
  $\WO^{2}_{d}(\Omega)$ we mean the subset of
$W^{2}_{d}(\Omega)$ consisting of continuous functions in $\bar\Omega$
vanishing on $\partial \Omega$.

In the parabolic part of the article 
we fix $T\in(0,\infty)$ and use the notation
$$
C=[0,T)\times B,
\quad \partial' C=\partial C
\setminus (\{0\}\times\bar B).
$$
We call $\partial' C$ the parabolic boundary of $C$.

Everywhere below
$\psi(x)=1-|x|^{2}$.

\mysection{Elliptic equations
of the main type in a ball}
                                      \label{section 8.11.1}

\begin{theorem}
                                             \label{theorem 8.6.1}

Let $u\in W^{2}_{d,\loc}(B)\cap C(\bar B)$ be a convex function in $B$,
 and let
 $\alpha\in[0,(d+1)/2)$. Then,
  for any $ x_{0}\in B $,
\begin{equation}
                                               \label{8.5.1}
  u(x_{0}) \geq \inf_{\partial B}u-N(d,\alpha)
\psi^{\beta}(x_{0})\Big(\int_{B}\psi^{\alpha} 
 \det D^{2}u  \,dx\Big)^{1/d},
\end{equation}
where
$$
\beta=(d+1-2\alpha)/(2d).
$$

\end{theorem}

\begin{corollary}
                                       \label{corollary 8.9.1}
Under the conditions of Theorem \ref{theorem 8.6.1},
if $u=0$ on $\partial B$, we have
\begin{equation}
                                               \label{8.9.5}
\sup_{B}|u|\leq N(\alpha,d)\Big(\int_{B}\psi^{\alpha} 
 \det D^{2}u \,dx\Big)^{1/d}.
\end{equation}
\end{corollary}

\begin{remark}
                                         \label{remark 8.5.1}
It might be that \eqref{8.5.1}  also holds if
$\alpha=(d+1)/2$. At least this is true indeed
if $d=1$. Generally,
 estimate \eqref{8.9.5} is close to be optimal in the following sense.
Take $\alpha>(d+1)/2$ and a sequence $x_{n}\in B$
such that $|x_{n}|\to1$ as $n\to\infty$. Then one can construct
a sequence of smooth in $\bar B$, convex
functions $u^{n}$, vanishing on the boundary,
such that $u^{n}(x_{n})\to-\infty$ and the
integral in the right-hand side
of \eqref{8.9.5} stays bounded.

Our argument showing this is rather  descriptive
dropping some rigorous justifications. But the author is sure
that the reader will be able to make it absolutely rigorous.
To construct such  a sequence of $u^{n}$,
define negative $v^{n}(x_{n})$ so that
$$
|v^{n}(x_{n})|^{d}=\psi^{(d+1)/2-\alpha}(x_{n})
$$
and introduce
 a cone with vertex at $ (x_{n},v^{n}(x_{n}))$
and base $ \partial B $. 
Let this cone be the graph of a function which we call
$v^{n}(x)$. Then mollify $v^{n}$ near $x_{n}$,
without changing it for $x$ not close to $x_{n}$
so that the new function, $u^{n}$, will be smooth, convex, and
close to $v^{n}$, so that 
$u^{n}(x_{n})\to-\infty$ (observe that 
$v^{n}(x_{n})\to-\infty$). 

Note that, since $u^{n}$ is smooth,
by change of variables formula, 
$$
\int_{B} 
 \det D^{2}u^{n} \,dx=|Du^{n}(B)|,
$$
where $Du^{n}(B)=\{Du^{n}(x):x\in B\}$.
It turns out that $ Du^{n}(B) $ and $|Du^{n}(B)|$ are
  independent of what we did with $v^{n} $
in the small neighborhood of $x_{n}$. Indeed,
if $x\in B$ and $p=Du^{n}(x)$, then the hyper-plane
$y=u^{n}(x)+p\cdot(z-x), z\in\bR^{d}$, can be shifted down,
if necessary,
so that it will become a supporting plane for the graph of $v^{n}$
at $(x_{n},v^{n}(x_{n}))$. Obviously, all supporting 
planes for the graph of $v^{n}$
at $(x_{n},v^{n}(x_{n}))$ can be obtained in this way,
so that $ Du^{n}(B)$ is just the collection of $p\in\bR^{d}$
such that $b+p\cdot(z-x_{n})$ is
a supporting plane for the graph of $v^{n}$
at $(x_{n},v^{n}(x_{n}))$ for a $b\in \bR$.

In analytic terms this means that
$$
p\in Du^{n}(B)\Longleftrightarrow
(p,x_{n})+b=v^{n}(x_{n}),\quad (p,x)+b\leq 0\quad\forall x\in\partial B.
$$
 We can rewrite the latter conditions
as
$$ 
{\rm max}_{|x|=1}(x-x_{n},p)\leq -v^{n}(x_{n}),\quad
|p|-(x_{n},p) 
\leq -v^{n}(x_{n}),
$$
$$ 
|p|^{2}-(x_{n},p)^{2}+2v^{n}(x_{n})(x_{n},p)\leq |v^{n}(x_{n})|^{2},
$$
 and $  b=v^{n}(x_{n})-(p,x_{n})$. 
For $ x_{n}=|x_{n}|e_{1}$ we have
$$ 
(1-|x_{n}|^{2})\bigg(p^{1}+
\frac{v^{n}(x_{n})|x_{n}|}{ 1-|x_{n}|^{2}} 
\bigg)^{2}
+\sum_{i\geq2}(p^{i})^{2}\leq 
\frac{|v^{n}(x_{n})|^{2}}{ 1-|x_{n}|^{2}}.
$$

It follows that $Du^{n}(B)$ is an ellipsoid
whose $ d-1$ principal semi-axes have length 
$$ 
|v^{n}(x_{n})|/\sqrt{1-|x_{n}|^{2}}
$$ 
and the remaining
one is of length 
$$ 
|v^{n}(x_{n})|/(1-|x_{n}|^{2}) .
$$ 
Hence,
$$
|Du^{n}(B)|=\omega_{d}|v^{n}(x_{n})|^{d}\psi^{-(d+1)/2}(x_{n}).
$$

After that observe that $\det D^{2}u^{n}\ne0$ only in the neighborhood
of $x_{n}$, where we changed $v^{n}$.
Then
$$
 \int_{B}\psi^{\alpha}
\det D^{2}u^{n} \,dx 
$$
  is close to 
$$
J:=\psi^{\alpha}(x_{n})\int_{B} 
\det D^{2}u^{n} \,dx=\psi^{\alpha}(x_{n})|Du^{n}(B)|
$$
$$
=\omega_{d}\psi^{\alpha-(d+1)/2}(x_{n})|v^{n}(x_{n})|^{d}=\omega_{d},
$$
and becomes as close   to $J$ as we wish as we shrink
the neighborhoods of $x_{n}$, where we changed $v^{n}$.
This finishes the proof of the claim
made at the beginning of this remark.  
 
\end{remark}
\begin{remark}
                                                    \label{remark 8.28.1}
It suffices to prove \eqref{8.5.1}
 with $B_{r}$
in place  of $B$ and $r^{2}-|x|^{2}$
in place of $\psi$ for $r<1$ close to 1. In that case
we have $u\in W^{2}_{d}(B_{r})$. This shows that in the proof of Theorem 
\ref{theorem 8.6.1}
without losing generality we may assume
that $u\in W^{2}_{d}(B)$.
\end{remark}

 \begin{lemma}
                                            \label{lemma 8.7.1}
Let $u\in W^{2}_{d}(B)$ be convex in $B$, $0\leq s<r\leq1$,
and suppose that $\det D^{2}u\ne0$ only on $B_{r}\setminus
B_{s}$ and $u=0$ on $\partial B$. 
  Then

(a) for $|x|\leq s$ we have
\begin{equation}
                                              \label{8.7.2}
|u(x)|\leq M\Big(\frac{1-s^{2}}{1-r^{2}}\Big)^{(d+1)/(2d)},
\end{equation}
where
$$
M=\omega_{d}^{-1/d}(1-r ^{2})^{(d+1)/(2d)}
\big(\int_{B_{r}\setminus B_{s}} 
\det D^{2}u \,dx\big)^{1/d},
$$

(b) for $|x|\geq r $ we have
\begin{equation}
                                              \label{8.7.3}
|u(x)|\leq M(1-|x|)/(1-r).
\end{equation}

\end{lemma}

Proof. (a) Since in $B_{s}$ we have
\begin{equation}
                                                 \label{8.17.1}
\inf_{a\in A}a^{ij}D_{ij}u=0,
\end{equation}
where $A$ is the set of $d\times d$ symmetric
nonnegative matrices with unit trace, by the maximum principle
$u$ in $\bar B_{s}$ attains its minimum on $\partial B_{s}$.
Furthermore,
the classical Aleksandrov estimate
(the derivation of which is just part of the arguments
in Remark \ref{remark 8.5.1}) says that, for any $x\in B$,
\begin{equation}
                               \label{9.5.1}
|u(x)|\leq \omega_{d}^{-1/d}(1-|x|^{2})^{(d+1)/(2d)}
\Big(\int_{B } 
 \det D^{2}u \,dx\Big)^{1/d},
\end{equation}
where the integral can be restricted to $B_{r}\setminus B_{s}$.
It follows that $|u(x)|$ on $\partial B_{s}$ is dominated by
the right-hand side of \eqref{8.7.2} and proves (a).

(b) Observe that in $C_{r}:=B\setminus\bar B_{r}$
the function $u$ satisfies \eqref{8.17.1}.
The function $v(x):=M(|x|-1)/(1-r)$
also satisfies this equation in $C_{r}$
and is less than $u$ on $\partial C_{r}$.
By Theorem 4.1.18 of \cite{Kr_18} we have $u\geq v$
in $C_{r}$ and this is \eqref{8.7.3}.
The lemma is proved.

Next, we need a special partition of unity
in $B$. For $n=0,1,2,...$ introduce
 $r_{n}=1-e^{-n}$, $r_{-1}=0$, and
find nonnegative $C^{\infty}$-functions
$\zeta_{n}$ on $[0,1]$ such that $\zeta_{n}\leq 1$,
$$
\zeta_{n}=1\quad\text{on} \quad [r_{n},r_{n+1}],
\quad \zeta_{n}=0\quad\text{for} \quad t< (r_{n-1}+r_{n})/2=:s_{n-1},
$$
$$
  \zeta_{n}=0\quad\text{for} \quad t\geq (r_{n+1}+r_{n+2})/2=:s_{n+1}.
$$
After that set
$$
\eta=\sum_{n=0}^{\infty}\zeta_{n},\quad\eta_{n}=\zeta_{n}/\eta.
$$
 Clearly, $\eta$ is infinitely differentiable on
$[0,1)$ and $0\leq\eta\leq 3$.

\begin{lemma}
                                                \label{lemma 8.8.1}
Let $u$ be four times continuously
differentiable in $\bar B$ function,
which is convex in $B$, is nonnegative on $\partial B$
and is such that $f:= \det D^{2}u >0$
in $\bar B$. Introduce $u_{n}$
as $C^{1,1}(\bar B)$-functions which are
convex in $B$, vanish  on $\partial B$ and
satisfy
$$
\det D^{2}u_{n}=f\eta_{n}^{d}
$$
in $B$ (a.e.). Then
$$
u\geq\sum_{n=1}^{\infty}u_{n}  .   
$$
\end{lemma}

Proof. First of all note that, since $g:=f^{1/d}$
is at least in $C^{1,1}(\bar D)$, the functions
$u_{n}$ with the described properties exist
and are unique according to   \S 4 in \cite{Kr_81}.

Recall that
 $u$ is a   smooth solution of
$$
\inf_{a\in A}[a^{ij}D_{ij}u-gd\sqrt[d]{ \det a}]=0
$$
nonnegative on $\partial B$, and $u_{n}$ are   unique 
$C^{1,1}$-solutions of
$$
\inf_{a\in A}[a^{ij}D_{ij}u_{n}- g\eta_{n}d\sqrt[d]{ \det a}]=0
$$
in $B$ (a.e.) vanishing on $\partial B$. Set $v_{n}=u_{0}+...+u_{n}$,
$g_{n}=g\eta_{0}+...+g\eta_{n}$ and note that
$$
\inf_{a\in A}[a^{ij}D_{ij}v_{n}- g_{n} d\sqrt[d]{ \det a}]\geq0
$$
in $B$ (a.e.). It follows that
$$
\sup_{a\in A}[a^{ij}D_{ij}(v_{n}-u)
- (g_{n}-g)d \sqrt[d]{ \det a}]\geq0
$$
and there is an $A$-valued function $a=a(x)$ such that
in $B$ (a.e.)
$$
a^{ij}D_{ij}(v_{n}-u)
- (g_{n}-g)d \sqrt[d]{ \det a} \geq0.
$$
By Theorem 3.3.4 of \cite{Kr_18} (see also the end of the introduction
into Chapter 3 in \cite{Kr_18})
$$
v_{n}-u\leq N(d)\|g_{n}-g\|_{L_{d}(B)}.
$$
This proves the lemma since the norm on the right
tends to zero as $n\to\infty$.

{\bf Proof of Theorem \ref{theorem 8.6.1}}.
According to Remark \ref{remark 8.28.1}
we may assume that $u\in W^{2}_{d}(B)$.
By having in mind approximations and adding to $u$ the function
$-\varepsilon\psi$ and then letting $\varepsilon
\downarrow0$, we convince ourselves that
we may also assume that $u$ is a strictly convex $C^{4}$-function. 
By replacing $u$ with $u-\inf_{\partial B}u$ we see that we may assume that $u\geq0$
on $\partial B$ as well.
Under these additional assumptions
fix $x_{0}\in B$  
and define $n_{0}$ as the
smallest $n\geq1$ such that $|x_{0}|\leq s_{n }$,
that is 
$$
|x_{0}|\leq s_{n_{0} },\quad  |x_{0}|\geq s_{n_{0}-1}.
$$
Below we are going to use a few times that
the ratio $(1-s_{n})/(1-s_{n+1})$
is bounded from above and away from zero by absolute
constants independent of $n$ and that, for
$|x|\in[s_{n-1},s_{n+2}]$, the ratio
$\psi^{\alpha}(x)/(1-s_{n})^{\alpha}$ 
is bounded from above and away from zero by  
constants independent of $n$ and depending only on $\alpha$.

Take $u_{n}$ from Lemma \ref{lemma 8.8.1} and observe that, by Lemma \ref{lemma 8.7.1}
(and \eqref{9.5.1}), 
if $n\leq n_{0}$, then
$$
|u_{n}(x_{0})|\leq N(d)(1-|x_{0}|)(1-s_{n})^{-1}
(1-s_{n}^{2})^{(d+1)/(2d)}
\Big(\int_{B } \eta_{n }^{d}
 \det D^{2}u  \,dx\Big)^{1/d}
$$
$$
\leq N(d)(1-|x_{0}|)(1-s_{n})^{(d+1)/(2d)-1-\alpha/d}
\Big(\int_{B } \eta_{n }^{d}\psi^{\alpha}
 \det D^{2}u \,dx\Big)^{1/d}
$$
$$
\leq N(d)(1-|x_{0}|)e^{n(1-\beta)}\Big(\int_{B } \eta_{n }^{d}\psi^{\alpha}
 \det D^{2}u  \,dx\Big)^{1/d}.
$$

By H\"older's inequality and in light of the fact that $1-\beta>0$ ($d\geq2$),
$$
\sum_{n=0}^{n_{0}}|u_{n}(x_{0})|
\leq N(d)(1-|x_{0}|)\Big(\sum_{n=0}^{n_{0}}e^{n(1-\beta) d/(d-1)}\Big)^{(d-1)/d}
$$
$$
\times \Big( \int_{B}\psi^{\alpha}
 \det D^{2}u \sum_{n=0}^{n_{0}}\eta_{n}^{d}\,dx\Big)^{1/d}
$$
$$
\leq N(d,\alpha)(1-|x_{0}|)e^{n_{0}(1-\beta)}
\Big( \int_{B}\psi^{\alpha}
 \det D^{2}u  \,dx\Big)^{1/d}
$$
\begin{equation}
                                                 \label{8.9.2}
\leq N(d,\alpha)\psi^{\beta}(x_{0})
\Big( \int_{B}\psi^{\alpha}
 \det D^{2}u  \,dx\Big)^{1/d}.
\end{equation}

If   $n>n_{0}$, then by Lemma \ref{lemma 8.7.1}
$$
|u_{n}(x_{0})|\leq N(d)(1-s_{n}^{2})^{(d+1)/(2d)}
\Big( \int_{B}\eta_{n}^{d}
 \det D^{2}u \,dx\Big)^{1/d}
$$
 $$
\leq N(d) e^{-n \beta}\big(\int_{B } \eta_{n }^{d}\psi^{\alpha}
\det D^{2}u \,dx\Big)^{1/d}.
$$
Since $\beta>0$,
$$
\sum_{n>n_{0}}|u_{n}(x_{0})|\leq N(d)\Big(\sum_{n>n_{0}}
e^{-n \beta d/(d-1)}\Big)^{(d-1)/d}
\Big(\int_{B}\psi^{\alpha}
\det D^{2}u \,dx\Big)^{1/d}
$$
$$
 \leq N(d,\alpha)
e^{-n_{0} \beta }\Big(\int_{B}\psi^{\alpha}
\det D^{2}u\,dx\Big)^{1/d}
$$
$$
\leq N(d,\alpha)\psi^{\beta}(x_{0})
\Big( \int_{B}\psi^{\alpha}
\det D^{2}u \,dx\Big)^{1/d}.
$$
Upon combining this with \eqref{8.9.2} and Lemma
\ref{lemma 8.8.1} we get \eqref{8.5.1}.
The theorem is proved.

 \begin{theorem} 
                        \label{theorem 8.11.1}

Let $ a(x)=(a^{ij}(x))$ be a $ d\times d$-symmetric
nonnegative definite matrix-valued
measurable function on $ B$
such that $ \tr a>0$ in $ B$. Let $\alpha\in[0,(d+1)/2)$
and $u\in W^{2}_{d,\loc}(B)\cap C(\bar B)$.
Introduce
$$ 
L_{0}u= a^{ij}D_{ij}u .
$$
Then, 
for any $ x_{0}\in B $, $(0/0:=0)$
\begin{equation}
                                                   \label{8.9.7}
u(x_{0}) \leq \sup_{\partial B}u+N(d,\alpha)\psi^{\beta}(x_{0})
\Big(\int_{B}\psi^{\alpha}I_{Lu<0}
(\det a)^{-1}|L_{0}u|^{d} \,dx\Big)^{1/d}.
\end{equation}
\end{theorem}

Proof. 
As in Remark \ref{remark 8.28.1}
we may assume that $u\in W^{2}_{d}(B)$.
Since $\tr a>0$, by the homogeneity of estimate
\eqref{8.9.7} we may assume that
\begin{equation}
                                              \label{8.11.3}
\tr a\equiv 1.
\end{equation}

{\em A particular case\/}. Suppose that $\det a\geq \varepsilon$,
where $\varepsilon>0$ is a constant. In that case it is easy
to pass to the limit in \eqref{8.9.7} from smooth $u$  to the ones in $W^{2}_{d}(B)$. Therefore,
we assume that $u\in C^{4}(\bar B)$. We may also assume
that $a$ is infinitely differentiable so that
$$
f:=-\frac{1}{d}(\det a)^{-1/d}L_{0}u
$$
is twice continuously differentiable in $\bar B$.
In that case, as we know, there exists a unique
$C^{1,1}(\bar B)$-function $v$ which is
concave  in $B$, vanishes  on $\partial B$, and
satisfies
\begin{equation}
                                                \label{8.8.1}
\det(- D^{2}v) =f_{+}^{d}
\end{equation}
in $B$ (a.e.). By the way, it is proved in \cite{GTW_99}
that the same result holds if 
we replace $f_{+}^{d}$ with $g^{d-1}$,
provided that $g\in C^{1,1}(\bar B), g\geq0$. Interestingly enough,
since  $f_{+}^{d/(d-1)}$ is not necessarily in $C^{1,1}$,
the  result in \cite{GTW_99} is not applicable here.

Then we know that $v$ also satisfies
\begin{equation}
                                                \label{8.11.5}
\sup_{a\in A}[a^{ij}D_{ij}v+fd\sqrt[d]{\det a}]=0
\end{equation}
and for $\bar u=u-\sup_{\partial B}u$ we obviously have
$$
\sup_{a\in A}[a^{ij}D_{ij}\bar u+fd\sqrt[d]{\det a}]\geq0.
$$
By Theorem 4.1.18 of \cite{Kr_18} we have $\bar u\leq v$
and a reference to Theorem \ref{theorem 8.6.1}
completes considering this case.

{\em General case\/}. In order  to drop the additional assumptions
 we take $\gamma,\delta>0$,
 take $\psi$  from Lemma 3.1.8 of \cite{Kr_18}
  introduce  $ L^{\delta}= L+\delta \Delta $,
and apply the  above result  to  
$u^{\gamma}=u-\gamma(\psi+1)$ and  $ L ^{\delta}$.
Then following almost  word for word  
the appropriate parts of the proof
of Lemma 3.2.4 of \cite{Kr_18}, we  arrive at  \eqref{8.9.7}. 
The theorem is proved.

 \mysection{General elliptic equations in $B$}
                                      \label{section 8.11.2}

Here we generalize Theorem \ref{theorem 8.11.1}
for equation with lower order terms.
\begin{theorem}
                         \label{theorem 8.11.2}
Under the assumptions of Theorem \ref{theorem 8.11.1}
let $b=(b^{i}(x))$ be $\bR^{d}$-valued
measurable
function on $B$  and $c=c(x)$ be a nonnegative bounded measurable
function on $B$ and let  $|b|\leq K\mu$ in $B$,
where $K$ is a fixed constant and $\mu=\mu(x)$
is the smallest eigenvalue of $a=a(x)$. Then
estimate \eqref{8.9.7} holds again with $N=N(d,\alpha,K)$
if we replace $L_{0}u$ in \eqref{8.9.7} with
\begin{equation}
                                              \label{8.11.2}
Lu=a^{ij}D_{ij}u+b^{i}D_{i}u-cu
\end{equation}
and $\sup_{\partial B}u$ with $\sup_{\partial B}u_{+}$.

\end{theorem}  

To prove the theorem we need a lemma.

\begin{lemma}
                                      \label{lemma 8.11.3}
Take $\gamma\in[0,1]$ and introduce
$$
w(x)=\int_{|x|}^{1}e^{Ks}\int_{0}^{s}\frac{e^{-Kt}}{
(1-t^{2})^{\gamma}}\,dtds.
$$
Then $w(x)$ has bounded first and second order derivatives in $B_{r}$
for any $r\in(0,1)$,
and satisfies there
\begin{equation}
                                              \label{8.11.4}
Lw+\frac{\mu}{(1-|x|^{2})^{\gamma}}\leq0.
\end{equation}
\end{lemma}

Proof. We have
$$
D_{i}w(x)=-\frac{x^{i}}{|x|} e^{K|x|}\int_{0}^{|x|}
\frac{e^{-Kt}}{
(1-t^{2})^{\gamma}}\,dt,
$$
$$
D_{ij}w(x)=-\frac{x^{i}x^{j}}{|x|^{2}}\frac{1}{(1-|x|^{2})^{\gamma}}
$$
$$
-\Big(\frac{\delta^{ij}}{|x|}-\frac{x^{i}x^{j}}{|x|^{3}}
+K\frac{x^{i}x^{j}}{|x|^{2}}\Big)e^{K|x|}
\int_{0}^{|x|}
\frac{e^{-Kt}}{
(1-t^{2})^{\gamma}}\,dt.
$$

Since $a^{ij}x^{i}x^{j}\geq \mu(x)|x|^{2}$ and
$a^{ij}x^{i}x^{j}\leq |x|^{2}a^{ij}\delta^{ij}$,
it follows that
$$
Lw\leq - \frac{\mu}{(1-|x|^{2})^{\gamma}}-cw
$$
$$
 -K\mu e^{K|x|}
\int_{0}^{|x|}
\frac{e^{-Kt}}{
(1-t^{2})^{\gamma}}\,dt
-\frac{b^{i}x^{i}}{|x|} e^{K|x|}\int_{0}^{|x|}
\frac{e^{-Kt}}{
(1-t^{2})^{\gamma}}\,dt.
$$
This easily implies \eqref{8.11.4}. The lemma is proved.

{\bf Proof of Theorem \ref{theorem 8.11.2}}. As in  the proof of 
Theorem \ref{theorem 8.11.1} we assume \eqref{8.11.3},
$u\in W^{2}_{d}(B)$, 
and first also assume that  $\det a\geq \varepsilon$,
where $\varepsilon>0$ is a constant. This allows us to
also assume that $u$, $a$, $b$, and $c$ are infinitely
differentiable so that
$$
f:=-\frac{1}{d}(\det a)^{-1/d}Lu
$$
is twice continuously differentiable in $\bar B$.

After that we define a function $v$ having the same properties
as   in  the proof of 
Theorem \ref{theorem 8.11.1} by solving
\eqref{8.8.1}. We also fix $x_{0}\in B$, take the function $w$ from Lemma
\ref{lemma 8.11.3} with $\gamma=1-\beta$, take $N(d,\alpha)$ from
\eqref{8.5.1} and set  
$$
N(d,\alpha)
\psi^{\beta}(x_{0})\Big(\int_{B}\psi^{\alpha} 
 \det D^{2}v  \,dx\Big)^{1/d}=:N_{1}\psi^{ \beta}(x_{0}).
$$ 

Observe that, as for any concave function vanishing on
$\partial B$, we have $|Dv|\leq v/(1-|x|)\leq 2v/\psi$.
Also note that for $\bar u=u-\sup_{\partial B}u_{+}$
we have
$
L\bar u\geq Lu$.
Then for the function
$$
\phi=v+2KN_{1}w
$$
owing to \eqref{8.5.1}, \eqref{8.11.4},
 and \eqref{8.11.5}, (also recall the definition of
$f$) we get
$$
L\bar u\geq Lv-b^{i}D_{i}v+cv\geq Lv-2|b| v/\psi
\geq Lv-2N_{1}|b| \psi^{-(1-\beta)}\geq L\phi.
$$
By the maximum principle $\bar u\leq\phi$ and, to finish the proof
in out particular case, it only remains
to observe that by using l'Hospital's rule
it is easy to check that
\begin{equation}
                                                \label{9.17.1}
\lim_{|x|\uparrow 1}\frac{w(x)}{(1-|x|^{2})^{ \beta}}
=2^{-\beta }  \beta  ^{-1}e^{K}\lim_{s\uparrow 1}  
(1-s)^{1-\beta}
\int_{0}^{s}\frac{e^{-Kt}}{
(1-t^{2})^{1-\beta}}\,dt=0.
\end{equation}

The general case is dealt with as at the end of the proof
of Theorem \ref{theorem 8.11.1}.
The theorem is proved.

\mysection{The case of general smooth domains}
                                      \label{section 8.15.1}

Let $\Omega$ be a $C^{1,1}$ bounded domain in $\bR^{d}$,
$a=(a^{ij}(x))$ be $d\times d$ symmetric matrix-valued
function on $\bR^{d}$, $b=(b^{i}(x))$ be $\bR^{d}$-valued
function on $\bR^{d}$, and $c=c(x)$ be a nonnegative bounded
function on $\bR^{d}$. We assume that these functions are measurable, fix three numbers $\delta,K >0$,
and assume that for all values of arguments and $\lambda\in\bR^{d}$
\begin{equation}
                                                        \label{8.15.5}
\delta^{-1}|\lambda|^{2}\geq
a^{ij}\lambda^{i}\lambda^{j}\geq \delta|\lambda|^{2},
\quad |b|\leq K, \quad  c\geq 0.
\end{equation}
Finally, take $\alpha\in[0,(d+1)/2)$ and define 
$\rho_{\Omega}(x)=\dist(x,\Omega^{c})$. 

\begin{theorem}
                                 \label{theorem 8.15.1}

Let $u\in W^{2}_{d,\loc}(\Omega)\cap C(\bar\Omega)$.
Then, 
for any $ x \in \Omega $, 
\begin{equation}
                                                   \label{8.15.1}
u(x ) \leq \sup_{\partial\Omega}u_{+}+N \rho_{\Omega}^{\beta}(x )
\|\rho_{\Omega}^{\alpha} 
  (Lu)_{-}\|_{L_{d}(\Omega)},
\end{equation}
where $L$ is taken from \eqref{8.11.2} and
the constant $N$ depends only on $\Omega$, $\delta$, $\alpha$, and $K$.
\end{theorem}

To prove this theorem we need three lemmas.
 
\begin{lemma}
                                             \label{lemma 8.11.1}
Set $B^{+}_{r}=B_{r}\cap\{x\in\bR^{d}:x^{1}>0\}$
and suppose that $u\in W^{2}_{d,\loc}(B^{+}_{1})\cap C(\bar
B^{+}_{1})$,
$u\leq0$ on $\partial B^{+}_{1}$. Also assume that
the  coefficients of $L$ are infinitely differentiable
and $f:=Lu \in L_{d}(B^{+}_{1})$.
 Then in $B^{+}_{1}$
we have
\begin{equation}
                                                        \label{8.15.4}
u(x ) \leq N (x^{1})^{\beta}
\|M^{\alpha} 
  f_{-}\|_{L_{d}(B^{+}_{1})},
\end{equation}
where $M^{\alpha}(x)=(x^{1})^{\alpha}$  and the constant
$N$ depends only on $d,\delta$, $\alpha$, and $K$.

\end{lemma}

Proof. Find
a $C^{1,1}$ domain $B'$ which contains $B^{+}_{2}$,
  has $B_{2}\cap\{x\in\bR^{d}:x^{1}=0\}$
as part of its boundary, and is $C^{1,1}$-diffeomorphic
to $B$. 

Then, extend $f$ as zero outside $B^{+}_{1}$
and denote by $v$ the 
function of class $\WO^{2}_{d}(B')$ satisfying
\begin{equation}
                                                        \label{8.15.6}
 Lv=-f_{-}
\end{equation}
in $B'$ (a.e.). It is a classical fact that
such a function exists and is unique.
By the maximum principle $v\geq 0$, $v\geq u$ in $B^{+}_{1}$.
Furthermore, is we apply the diffeomorphism
mentioned above to \eqref{8.15.6} then we will see
that the image $v'$ of $v$ will satisfy the equation
$$
 L'v'=-f'_{-}
$$
in $B$ (a.e.),
where $f'_{-}$ is the image of $f_{-}$ and $ L'$
is the image of $L$.

By Theorem \ref{theorem 8.11.2}
\begin{equation}
                                                        \label{8.15.7}
v'(x')\leq N\psi^{\beta}(x')
\|\psi^{\alpha/d} f'_{-}\|_{L_{d}(B)},
\end{equation}
where $\psi(x')$ is the distance from $x'$ to the 
boundary of $B$. If $x'$ is the image of $x$,
$\psi(x')$ is comparable to the distance of $x$ to
the boundary of $B'$, since the diffeomorphism
and its inverse are   Lipschitz continuous.

It only remains to write down
\eqref{8.15.7} in the original coordinates and
use the fact that $v\geq u$  and that for $x\in B^{+}_{1}$
it distance  to
the boundary of $B'$ equals $x^{1}$.
The lemma is proved.

Next we use a well-known fact
(see, for instance Lemma 8.8
in \cite{Kr_11})
that  there exists a function $\Psi\in C^{1,1}(\bar \Omega)$
such that, for a constant $N$ depending only on
$\Omega$, $\delta$, and $K$, we have on $\Omega$
\begin{equation}
                                                        \label{8.15.2}
 N\rho_{\Omega}\geq \Psi, \quad \Psi\geq \rho_{\Omega},
\end{equation}
\begin{equation}
                                                        \label{8.15.3}
L\Psi+c\Psi\leq-1.
\end{equation}

\begin{lemma}
                                            \label{lemma 8.15.1}
Let the coefficients of $L$ be infinitely differentiable.
Take $\rho_{0}\in(0,1)$, $x_{0}\in \Omega$ and suppose that
$2\rho:=\rho_{\Omega}(x_{0})\geq  \rho_{0}$. Let $\gamma\in (0,1]$ and let
$\Phi$ be the classical
 solution of $L\Phi=0$ in $B_{\rho }(x_{0})$
with boundary condition $\Phi=\Psi^{\gamma}$ on 
$\partial B_{\rho }(x_{0})$. Then
$$
\Phi(x_{0})\leq[1-\varepsilon \rho_{0}]^{\gamma} \Psi^{\gamma}(x_{0}),
$$
where $\varepsilon>0$ depends only on $K$, $\delta$, and $\Omega$.
\end{lemma}
 
Proof. A simple argument based on the maximum principle shows
that it suffices to concentrate on $\gamma=1$. 
In that case first we note that 
by the maximum principle $\Phi\geq0$. Therefore, $(L+c)\Phi\geq0$
and
$v:=\Psi-\Phi$ satisfies $(L+c)v\leq-1$ and $v=0$
on 
$\partial B_{\rho }(x_{0})$. Then elementary barriers 
show that $v(x_{0})\geq \varepsilon \rho ^{2}$,
where $\varepsilon>0$ depends only on $K$, $\delta$,
and the diameter of $\Omega$.
Since $\rho_{\Omega}(x_{0})$ and $\Psi(x_{0})$
are comparable we have (we use $\varepsilon$ as a generic
constant $>0$ depending only on $K$, $\delta$, and $\Omega$)
$$
\rho^{2}=\rho(1/2)\rho_{\Omega}(x_{0})\geq
\varepsilon \rho_{0}\Psi(x_{0}),
$$
$$
\Psi(x_{0})-\Phi(x_{0})\geq \varepsilon \rho_{0}\Psi(x_{0}),
$$
and the lemma is proved.

Below by $\rho_{0}$ we mean a number $>0$
such that  any point $\bar x_{0}\in\partial\Omega$
is the only common point of $\partial\Omega$
and the closure of a ball, say $B_{\rho_{0}}(y_{0})$,
 belonging to $\Omega$
with radius $\rho_{0}$.
Since $\Omega\in C^{1,1}$, such $\rho_{0}>0$ exists.
We further decrease $\rho_{0}$, if necessary, so that
there is a $C^{1,1}$-diffeomorphism with its first- and second-order derivatives
and the first- and second-order derivatives of its inverse
bounded by a constant depending only on $\Omega$ and
mapping
$B_{2\rho_{0}}(\bar x_{0})\cap \Omega$ onto
$B^{+}_{1}$ and $\bar B_{2\rho_{0}}(\bar x_{0})\cap 
\partial\Omega$ onto $\bar B^{+}_{1}\cap\{x^{1}=0\}$.

In the following lemma we consider the case that
$ \rho_{\Omega}(x_{0})\leq  \rho_{0}$ and denote
by $\bar x_{0}$ a   point on $\partial \Omega$
such that $\rho_{\Omega}(x_{0})=|\bar x_{0}-x_{0}|$.
By assumption there exists $y_{0}\in\Omega$ such that
$$
\bar B_{\rho_{0}}(y_{0})\cap\partial\Omega=\{\bar x_{0}\},
$$
and since both balls $B_{\rho_{0}}(y_{0})$ and
$B_{\rho_{\Omega}(x_{0})}(x_{0})$ touch $\partial\Omega$
at $\bar x_{0}$ and the former ball has a smaller radius,
the points $\bar x_{0},x_{0}$, and $y_{0}$ lie
on the same line and
\begin{equation}
                                                \label{8.16.1}
\dist(x_{0},B^{c}_{\rho_{0}}(y_{0}))=\rho_{\Omega}(x_{0}).
\end{equation}

\begin{lemma}
                                            \label{lemma 8.15.2}
Let the coefficients of $L$ be infinitely differentiable.
Take $\rho_{0}\in(0,1)$, $x_{0}\in \Omega$ and suppose that
$\rho:= \rho_{\Omega}(x_{0})\leq  \rho_{0}$. 
Take $y_{0}$ introduced before the lemma.
Let $\gamma\in (0,1]$ and let
$\Phi$ be the classical
 solution of $L\Phi=0$ in $B_{ \rho_{0} }(y_{0})$
with boundary condition $\Phi=\Psi^{\gamma}$ on 
$\partial B_{ \rho_{0} }(y_{0})$. Then
$$
\Phi(x_{0})\leq[1-\varepsilon \rho_{0}]^{\gamma} \Psi^{\gamma}(x_{0}),
$$
where $\varepsilon>0$ depends only on $K$, $\delta$, and $\Omega$.
\end{lemma}

Proof. As in the proof of Lemma \ref{lemma 8.15.1},
we concentrate on the  case that $\gamma=1$ and
we have $(L+c)v\leq-1$. Simple barriers show that
(recall that $\rho_{0}$ is fixed)
$v(x)\geq \varepsilon(\rho_{0}-|x-y_{0}|)$ in $B_{ \rho_{0} }(y_{0})$,
where $\varepsilon>0$ depends only on $d$, $\delta$, $K$, and
$\rho_{0}$. Owing to \eqref{8.16.1} this shows that
$v(x_{0})\geq \varepsilon \rho_{\Omega}(x_{0})$
and we are done because $\rho_{\Omega}(x_{0})$ and $\Psi(x_{0})$
are comparable. The lemma is proved.

{\bf Proof of Theorem \ref{theorem 8.15.1}}.
An argument similar to the one in the beginning
of the proof of Theorem \ref{theorem 8.11.1}
shows that without losing generality we may assume 
that $u\in W^{2}_{d}(\Omega)$. After that
we can certainly concentrate on the case that
$u$ and the coefficients of $L$ are infinitely differentiable
in $\bar \Omega$. 
Then define $v$ as a unique classical solution of $Lv=Lu$ in $\Omega$ with zero boundary condition. 
By the maximum  principle $u-v\leq \sup_{\partial\Omega}u_{+}$. It follows that we only need
to estimate $v$ and consequently we may assume that
$u=0$ on $\partial\Omega$.

In that case the function
$$
v=\frac{u}{\Psi^{\beta}}
$$
is continuous in $\bar\Omega$, equals zero on $\partial\Omega$ and, hence,
attains its maximum value at a point $x_{0}\in \Omega$:
\begin{equation}
                       \label{8.16.3}
M:=\frac{u(x_{0})}{\Psi^{\beta}(x_{0})}\geq \frac{u(x)}{\Psi^{\beta}(x)}
\quad\forall x\in \Omega.
\end{equation}
If $M$ is less than zero, we have nothing
to prove. Therefore, we assume that
$$
u(x_{0})>0
$$ 
and consider two cases:

(a) $\rho_{\Omega}(x_{0})\geq \rho_{0}$,

(b) $\rho_{\Omega}(x_{0})< \rho_{0}$.

In case (a), in $B':=B_{\rho_{\Omega}(x_{0})/2}(x_{0})$ we have $u=v+h$,
where $v$ is the classical solution of $Lv=Lu$ 
in $B'$ with zero boundary value and $h$ is the classical
solution of $Lh=0$ in $B'$  with boundary condition
$h=u$ on $\partial B'$. We apply the Aleksandrov
estimate to $v$ and take into account that
on $B'$, $\rho_{\Omega}$ and $\Psi$ are comparable to a constant one.
Then we see that
\begin{equation}
                                                      \label{8.16.2}
v(x_{0})\leq N\| 
  (Lu)_{-}\|_{L_{d}(B')}\leq
N \Psi^{\beta}(x_{0} )
\|\rho_{\Omega}^{\alpha} 
  (Lu)_{-}\|_{L_{d}(\Omega)}.
\end{equation}

In what concerns $h$, observe that owing to \eqref{8.16.3}, 
by the maximum principle,
 it is less
than the solution $w$ of the equation $Lw=0$ in $B'$
with boundary condition $M\Psi^{\beta}$. By Lemma \ref{lemma 8.15.1},
$h(x_{0})\leq w(x_{0})\leq M\varepsilon\Psi^{\beta}(x_{0})$,
where $\varepsilon\in(0,1)$ depends only on $K,\delta,\alpha$,
and $\Omega$. It follows that
$$
M\leq N\|\rho_{\Omega}^{\alpha} 
  (Lu)_{-}\|_{L_{d}(\Omega)}+M\varepsilon,\quad
M\leq N\|\rho_{\Omega}^{\alpha} 
  (Lu)_{-}\|_{L_{d}(\Omega)},
$$
which yields \eqref{8.15.1}.

In case (b), in the ball $B_{\rho_{0}}(y_{0})$
introduced before Lemma \ref{lemma 8.15.2} we have
$u=v+h$, where $v$ is the solution of $Lv=Lu$
in $B_{\rho_{0}}(y_{0})$ with zero boundary condition
and $h$ satisfies $Lh=0$ and equals $u$ on
$\partial B_{\rho_{0}}(y_{0})$. As in case (a)
by the maximum principle and Lemma \ref{lemma 8.15.2}
we have $h(x_{0})\leq M\varepsilon \Psi^{\beta}(x_{0})$.

In what concerns $v$ observe that in $B_{\rho_{0}}(y_{0})$
 by the maximum
principle it is less than $w$ defined as
$W^{2}_{d,\loc}(B')\cap C(\bar B')$-solution of
$Lw=-f_{-}$ in
$$
B':=B_{2\rho_{0}}(\bar x_{0})\cap \Omega
$$   
vanishing on its boundary, where
$$
f:=I_{B_{\rho_{0}}(y_{0})}Lu.
$$
By an argument similar to the one used in the proof of Lemma
\ref{lemma 8.11.1} we obtain that
$$
w(x_{0})\leq N\rho^{\beta}_{\Omega}(x_{0})\|\rho_{\Omega}^{\alpha}f_{-}
\|_{L_{d}(B')}
\leq N\rho^{\beta}_{\Omega}(x_{0})\|\rho_{\Omega}^{\alpha}f_{-}
\|_{L_{d}(\Omega)}.
$$
Since $\rho_{\Omega}$ and $\Psi$ are comparable we conclude
$$
M\leq M\varepsilon+N\|\rho_{\Omega}^{\alpha}f_{-}
\|_{L_{d}(\Omega)}
$$
and this proves the theorem.

\mysection{Estimates for stochastic integrals}
                                                   \label{section 9.1.1}

Let $(\Omega,\cF,P)$ be a complete probability
space, $\{\cF_{t},t\geq0\}$ be an increasing filtration
of $\sigma$-fields $\cF_{t}\subset\cF$ each of 
which is complete relative to $\cF,P$. Suppose
that on $(\Omega,\cF,P)$ we are given a $d_{1}$-dimensional Wiener process which is $\cF_{t}$-adapted
and such that $w_{t}-w_{s}$ are independent of $\cF_{s}$
as long as $0\leq s\leq t<\infty$.

Fix some constants $\delta\in(0,1]$ and $K\geq0$.
Let $\sigma_{t}=\sigma_{t}(\omega)$ be a progressively
measurable with respect to $\{\cF_{t},t\geq0\}$,
$d\times d_{1}$-matrix valued process such that
$$
\delta^{-1}|\lambda|^{2}\geq a_{t}^{ij}\lambda^{i}\lambda^{j}
\geq\delta|\lambda|^{2}
$$
for all $\lambda\in\bR^{d}$, $t\geq0$, and
$\omega\in\Omega$, where $a=(1/2)\sigma\sigma^{*}$.

Let $b_{t} $ be a progressively
measurable with respect to $\{\cF_{t},t\geq0\}$,
$\bR^{d}$-valued process such that
$|b_{t}|\leq K$
for all  $t\geq0$ and
$\omega\in\Omega$.

Introduce
$$
x_{t}=\int_{0}^{t}\sigma_{s}\,dw_{s}+
\int_{0}^{s}b_{s}\,ds.
$$

Let $G$ be a bounded $C^{1,1}$ domain in $\bR^{d}$
containing the origin, set
$$
\rho_{G}(x)=\dist(x,G^{c}),
$$
take a nonnegative Borel function $f$ on $G$,
a number
$$
\alpha\in[0,(d+1)/2),
$$
and set
$$
\beta=(d+1-2\alpha)/(2d).
$$
Here is the result of this section.

\begin{theorem}
                                                  \label{theorem 9.1.1}
There exists a constant $N$, depending only on $G$, $\delta$, $\alpha$, and $K$, such that
\begin{equation}
                     \label{9.1.2}
E\int_{0}^{\tau}f(x_{t})\,dt
\leq N\rho_{G}^{\beta}(0)
\Big(\int_{G}\rho_{G}^{\alpha}f^{d}(x)\,dx\Big)^{1/d},
\end{equation}
where $\tau$ is the first exit time of $x_{t}$ from $G$.
\end{theorem}

Proof. A usual measure-theoretic argument
shows that it suffices to prove \eqref{9.1.2}
for bounded $f$.
Let $A$ be the collection of couples
$(a,b)$, where $a$ is a symmetric $d\times d$
matrix and $b\in \bR^{d}$ such that
$$
|b|\leq K,\quad \delta^{-1}|\lambda|^{2}\geq a ^{ij}\lambda^{i}\lambda^{j}
\geq\delta|\lambda|^{2}
$$
for all $\lambda\in\bR^{d}$.

As it follows from \cite{Wi09} or \cite{Kr_18}
there exists a unique $u\in W^{2}_{d}(G)$
vanishing on $\partial G$ and satisfying
$$
\sup_{(a,b)\in A}[a^{ij}D_{ij}u+b^{i}D_{i}u+f]=0
$$
in $G$ (a.e.). By It\^o's formula
(see, for instance, \cite{Kr_77})
$$
u(0)=E\int_{0}^{\tau}[-a^{ij}_{t}D_{ij}u(x_{t})-b^{i}_{t}
D_{i}u(x_{t})]\,dt\geq
E\int_{0}^{\tau}f(x_{t})\,dt.
$$

After that it only remains to 
apply Theorem \ref{theorem 8.15.1} first
observing that
there exist a measurable $A$-valued function
 $(a(x),b(x))$ such that
$$
a^{ij}(x)D_{ij}u(x)+b^{i}(x)D_{i}u(x)+f(x)=0
$$
in $G$ (a.e.).  The theorem is proved.

If Borel $\Gamma\subset G$, then
$$
G(\Gamma)=E\int_{0}^{\tau}I_{\Gamma}(x_{t})\,dt
$$
(the so-called Green's measure) is the mean time that the process
$x_{t},t\in[0,\tau]$, spends in $\Gamma$,
the mean time it occupies $\Gamma$
before exiting from $G$.
By taking $f=I_{\Gamma}$ in \eqref{9.1.2}
we come to the following.

\begin{corollary}
We have
$$
G(\Gamma)\leq 
N\rho_{G}^{\beta}(0)
\Big(\int_{\Gamma}\rho_{G}^{\alpha} (x)\,dx\Big)^{1/d}.
$$
\end{corollary}
\begin{remark}
By analyzing the arguments in Remark \ref{remark 8.5.1} and the proof of Theorem \ref{theorem 8.11.1}, it is not hard to show that in the whole
class of processes satisfying the conditions of Theorem \ref{theorem 9.1.1} with arbitrary
$\delta$ and $K=0$ estimate \eqref{9.1.2}
fails to hold if one replaces $d$ with $p<d$
even if $\alpha=0$ and also fails to hold for
  $\alpha>(d+1)/2$.
\end{remark}

\mysection{Parabolic equations
of the main type in a round cylinder}
                                      \label{section 8.11.10}

If $\sfu$ is a $d\times d$ symmetric matrix,
by $A_{ij}[\sfu]$ we denote the co-factor of $\sfu_{ij}$
in $\det\sfu$.

\begin{lemma}
                                            \label{lemma .8.18.1}
Let $u,v\in W^{1,2}_{d+1}(C)$, $u=v=0$ on $\partial'C$,
$u$ and $v$ be convex with respect to $x$  and satisfying $u\geq v$ 
in $C$. Then
\begin{equation}
                                                   \label{8.18.3}
\int_{C}\partial_{t}u\det D^{2}u\,dxdt\leq
 \int_{C}\partial_{t}v\det D^{2}v\,dxdt.
\end{equation}

\end{lemma}

Proof. By having in mind approximations, we may assume that $u$ and $v$ are infinitely
differentiable in $\bar C$. Then define $v_{\tau}=(1-\tau)v+\tau u$
and observe that, to prove \eqref{8.18.3}, it suffices to prove that
\begin{equation}
                                                   \label{8.18.4}
  \frac{d}{d\tau}\int_{C}\partial_{t}v_{\tau}\det D^{2}v_{\tau}\,dxdt\leq
 0.
\end{equation}

By denoting $\dot v_{\tau}=d v_{\tau}/d\tau$  ($u-v\geq0$), we see
that the left-hand side of \eqref{8.18.4} equals
$$
\int_{C}\partial_{t}\dot v_{\tau}\det D^{2}v_{\tau}\,dxdt
+\int_{C}\partial_{t}v_{\tau}D_{ij}\dot v_{\tau}A_{ij}[D^{2}v_{\tau}]\,dxdt=:I_{1}+I_{2},
$$
   By taking into account that
$$
\sum_{i}\frac{\partial}{\partial x^{i}}A_{ij}[D^{2}u]
=\sum_{i}\frac{\partial}{\partial x^{i}}A_{ji}[D^{2}u]=0
$$
for any $j$ and smooth $u$ and that $\partial_{t}v_{\tau}(t,x)=0$ for $|x|=1$,
we conclude that
$$
I_{2}=\int_{C}\partial_{t}D_{ij}v_{\tau}\dot v_{\tau}A_{ij}[D^{2}u]\,dxdt
=\int_{C}\dot v_{\tau}\partial_{t}\det D^{2}v_{\tau}\,dxdt.
$$

It follows that the left-hand side of \eqref{8.18.4} equals
$$
\int_{C}\partial_{t}\Big[\dot v_{\tau}\det D^{2}v_{\tau}\Big]\,dxdt
=-\int_{B}\dot v_{\tau}\det D^{2}v_{\tau}(0,x)\,dx.
$$
Since $\dot v_{\tau}\geq0$  and $\det D^{2}v_{\tau}\geq0$ we come to 
\eqref{8.18.4} thus proving the lemma.

\begin{theorem}
                                             \label{theorem 8.24.1}
Let $v \in W^{1,2}_{d+1}(C)$, $v =0$ on $\partial'C$, and
$v$   be convex with respect to $x$. Then,
  for any $ x_{0}\in B $,
\begin{equation}
                                                   \label{8.24.2}
|v(0,x_{0})|^{d+1}\leq (d+1)\omega_{d}^{-1} (1-|x_{0}|^{2})^{(d+1)/2}
 \int_{C}\partial_{t}v\det D^{2}v\,dxdt,
\end{equation}
where $\omega_{d}$ is the volume of $B$.

\end{theorem}

Proof. We may assume that $v$ is smooth.
 Define  $u (x_{0})=-1$ 
and introduce
 a cone with vertex at $ (x_{0},u(x_{0}))$
and base $ \partial B $. 
Let this cone be the graph of a function which we call
$u(x)$. 
Then mollify $u$
in $x$ near $x_{0}$ so that it remains convex, becomes smooth,
 larger than $u$ and coincides with $u$ in $B$ apart from a small
neighborhood of $x_{0}$. Call the resulting function $w$ and observe that $w(x_{0})$
is as close to $-1$ as we wish.

Then set
$$
w(t,x)=|v(t,x_{0})|w(x).
$$
By convexity of $v$ we have $w(t,x)\geq |v(t,x_{0})|u(x)\geq v(t,x)$. 
Also 
$$
\partial_{t}w\det D^{2}w(t,x)=-|v(t,x_{0})|^{d}\partial_{t}v(t,x_{0})
|w(x)|
\det D^{2}w(x).
$$
It follows that
$$
\int_{C}\partial_{t}w\det D^{2}w\,dxdt=
\frac{1}{d+1}|v(0,x_{0})|^{d+1}\int_{B}
|w(x)|
\det D^{2}w(x)\,dx,
$$
where  the integral on the right can be restricted to  the small
neighborhood of $x_{0}$ where we modified $u$ and hence this 
integral is as close as we wish to
$$
\int_{B}
\det D^{2}w(x)\,dx.
$$
The value of the last integral is well known to be
$\omega_{d} (1-|x_{0}|^{2})^{-(d+1)/2}$
(see, for instance, Remark \ref{remark 8.5.1})

After that it only remains to remember that
 by construction and Lemma~\ref{lemma .8.18.1}
$$
\int_{C}\partial_{t}w\det D^{2}w\,dxdt\leq
 \int_{C}\partial_{t}v\det D^{2}v\,dxdt.
$$
The theorem is proved.

Recall that $\psi(x)=1-|x|^{2}$.

\begin{theorem}
                                             \label{theorem 8.6.10}
Let $u \in W^{1,2}_{d+1,\loc}(C)\cap C(\bar C)$, 
$u$   be convex with respect to $x$ and increasing in $t$.
Let 
 $\alpha\in[0,(d+1)/2)$. Then,
  for any $ x_{0}\in B $,
\begin{equation}
                                               \label{8.5.10}
u(0,x_{0})\geq \inf_{\partial'C}u-N(d,\alpha)
\psi^{\beta}(x_{0})\Big(\int_{C}\psi^{\alpha} 
\partial_{t}u \det D^{2}u  \,dxdt\Big)^{1/(d+1)},
\end{equation}
where
$$
\beta=\frac{1}{2}-\frac{\alpha}{d+1}.
$$
\end{theorem}

\begin{corollary}
                                       \label{corollary 8.9.10}
Under the conditions of Theorem \ref{theorem 8.6.10},
if $u=0$ on $\partial' C$, we have
\begin{equation}
                                               \label{8.9.50}
\sup_{C}|u|\leq N(d,\alpha)\Big(\int_{C}\psi^{\alpha} 
\partial_{t}u \det D^{2}u\,dxdt\Big)^{1/(d+1)}.
\end{equation}
\end{corollary}

 We prove Theorem \ref{theorem 8.6.10} after some preparations.
For $1\geq r  > 0$ denote
$$
B_{r}=\{x:|x|< r\},\quad C_{r}=[0,T)\times B_{r}.
$$

 \begin{lemma}
                                            \label{lemma 8.7.10}
Let $v \in W^{1,2}_{d+1}(C)$, $v =0$ on $\partial'C$,  
$v$   be convex with respect to $x$ and increasing in $t$.
Let $s\in[0,r)$ and suppose that $\partial_{t}v\det D^{2}v\ne0$ only on $C_{r}\setminus
C_{s}$. 
  Then

(a) for $t\geq0,|x|\leq s$ we have
\begin{equation}
                                              \label{8.7.20}
|v(t,x)|\leq M\Big(\frac{1-s^{2}}{1-r^{2}}\Big)^{1/2},
\end{equation}
where
$$
M=\omega_{d}^{-1/(d+1)}(1-r ^{2})^{1/2}
\Big(\int_{C_{r}\setminus C_{s}} 
\partial_{t}v \det D^{2}v \,dxdt\Big)^{1/(d+1)},
$$

(b) for $|x|\geq r $ we have
\begin{equation}
                                              \label{8.7.30}
|v(t,x)|\leq M(1-|x|)/(1-r).
\end{equation}

\end{lemma}

Proof. (a) Observe that in $B_{s}$ we have
\begin{equation}
                                                 \label{8.17.10}
\inf_{(r,a)\in A}[r\partial_{t}v+a^{ij}D_{ij}v]=0,
\end{equation}
where $A$ is the set of couples $(r,a)$, where $a$ are $d\times d$ symmetric
nonnegative matrices, $r\geq0$, and $r+\tr a=1$. By the maximum principle
$u$ in $\bar C_{s}$ attains its minimum  on $[0,\infty)\times
\partial B_{s}$.
Furthermore, by Theorem \ref{theorem 8.24.1}
for any, $t_{0}\geq0,x_{0}\in B$, the quantity 
$|v(t_{0},x_{0})|^{d+1}$ is less than the right-hand side
of \eqref{8.24.2}
where the integral can be restricted to $C_{r}\setminus C{s}$.
It follows that $|v(t_{0},x_{0})|$ on $[0,\infty)\times\partial B_{s}$ is dominated by
the right-hand side of \eqref{8.7.20} and proves (a).

(b) Observe that in $G_{r}:=C\setminus\bar C_{r}$
the function $v$ satisfies \eqref{8.17.10}.
The function $w(t,x):=M(|x|-1)/(1-r)$
also satisfies this equation in $C_{r}$
and is less than $v(t,x)$     if $|x|=r$
or $|x|=1$. Also $v(T,x)=0\geq w(x)$. 
By Theorem 4.1.11 of \cite{Kr_18} we have $u\geq w$
in $G_{r}$ and this is \eqref{8.7.30}.
The lemma is proved.

Next, we use the special partition of unity
in $B$ introduced before Lemma \ref{lemma 8.8.1}. 

\begin{lemma}
                                                \label{lemma 8.8.10}
Let $u \in W^{1,2}_{d+1}(C)$,   
$u$   be convex with respect to $x$ and increasing in $t$.
Assume that
$$
\partial_{t}u\det D^{2}u=f_{+}^{d+1},
$$
in $C$ (a.e.), where $f\in L_{d+1}(C)\cap W^{1,2}_{\infty}(C_{r})$ for any $r\in(0,1)$,
and $f(T,x)=0$ in $B$.
  Introduce $u_{n}$
as $W^{1,2}_{\infty}(C)$-functions which are
convex in $x$,  increase in $t$, vanish 
on $\partial'C$, and
satisfy
$$
\partial_{t}u_{n}\det D^{2}u_{n}=f_{+}^{d+1}\eta_{n}^{d+1}
$$
in $C$ (a.e.). Then
$$
u\geq \inf_{\partial'C}u+\sum_{n=1}^{\infty}u_{n}    .  
$$
\end{lemma}

Proof. First of all note that, since $\psi^{-1/(d+1)}f \eta_{n}$
is   in $W^{1,2}_{\infty}( C)$ and vanish for $t=T$, the functions
$u_{n}$ with the described properties exist
and are unique according to   \S 5 in \cite{Kr_81}.

Recall that
 $\bar u:=u-\inf_{\partial'C}u$ satisfies
$$
\inf_{(r,a)\in A}[r\partial_{t}\bar u+a^{ij}D_{ij}\bar u-f(d+1)\sqrt[d+1]{r \det a}]=0
$$
 in $C$,  and $u_{n}$ are   unique 
$W^{1,2}_{\infty}(C)$-solutions of
$$
\inf_{(r,a)\in A}[r\partial_{t}u_{n}+a^{ij}D_{ij}u_{n}- f\eta_{n}
(d+1)\sqrt[d+1]{r \det a}]=0
$$
in $C$   vanishing on $\partial'C$. Set $v_{n}=u_{0}+...+u_{n}$,
$f_{n}=f\eta_{0}+...+f\eta_{n}$ and note that
$$
\inf_{(r,a)\in A}[r\partial_{t}v_{n}+a^{ij}D_{ij}v_{n}- f_{n}
 (d+1)\sqrt[d+1]{ r\det a}]\geq0
$$
in $C$. It follows that
$$
\sup_{(r,a)\in A}[r\partial_{t}(v_{n}-\bar u)+a^{ij}D_{ij}(v_{n}-\bar u)
- (f_{n}-f)(d+1) \sqrt[d+1]{r \det a}]\geq0
$$
and there is an $A$-valued function $(r,a)=(r(x),a(x))$ such that
in $C$ (a.e.)
$$
r\partial_{t}(v_{n}-\bar u)+a^{ij}D_{ij}(v_{n}-\bar u)
- (f_{n}-f)(d+1) \sqrt[d+1]{ r\det a} \geq0.
$$
By Theorem 3.2.3 of \cite{Kr_18} (see also the end of the introduction
into Chapter 3 in \cite{Kr_18})
$$
v_{n}-\bar u\leq N(d)\|f_{n}-f\|_{L_{d+1}(C)}.
$$
This proves the lemma since the norm on the right
tends to zero as $n\to\infty$.

\begin{lemma}
                                \label{lemma 8.25.1}
Let $u \in W^{1,2}_{d+1}(C)$,  
$u$   be convex with respect to $x$ and increasing in $t$.
Assume that
$$
\partial_{t}u\det D^{2}u=f_{+}^{d+1} 
$$
in $C$ (a.e.), where $f\in L_{d+1}(C)\cap  W^{1,2}_{\infty}(C_{r})$ for any $r\in(0,1)$,
and $f(T,x)=0$ in $B$. Let $\alpha\in[0,(d+1)/2)$. Then,
  for any $ x_{0}\in B $, estimate \eqref{8.5.10} holds. 
\end{lemma}

Proof. 
Fix $x_{0}\in B$  
and define $n_{0}$ as the
smallest $n\geq1$ such that $|x_{0}|\leq s_{n }$,
that is 
\begin{equation}
                                         \label{8.8.4}
|x_{0}|\leq s_{n_{0} },\quad  |x_{0}|\geq s_{n_{0}-1}.
\end{equation}
Below we are going to use a few times that
the ratio $(1-s_{n})/(1-s_{n+1})$
is bounded from above and away from zero by absolute
constants independent of $n$ and that, for
$x\in[s_{n-1},s_{n+2}]$, the ratio
$\psi^{\alpha}(x)/(1-s_{n})^{\alpha}$ 
is bounded from above and away from zero by  
constants independent of $n$ and depending only on $d$.

Take $u_{n}$ from Lemma \ref{lemma 8.8.10}
and observe that, by Lemma \ref{lemma 8.7.10}
(and Theorem \ref{theorem 8.24.1}), 
if $n\leq n_{0}$, then
$$
|u_{n}(0,x_{0})|\leq N(d)(1-|x_{0}|)(1-s_{n})^{-1}
(1-s_{n}^{2})^{1/2}
\Big(\int_{C} \eta_{n }^{d+1}
f_{+}^{d+1} \,dxdt\Big)^{1/(d+1)}
$$
$$
\leq N(d)(1-|x_{0}|)(1-s_{n})^{-1/2-\alpha/(d+1)}
\Big(\int_{C}  \psi^{\alpha}
\eta_{n }^{d+1}
f_{+}^{d+1}\,dxdt\Big)^{1/(d+1)}
$$
$$
\leq N(d)(1-|x_{0}|)e^{n(1-\beta)}\Big(\int_{C} 
 \psi^{\alpha}
\eta_{n }^{d+1}
f_{+}^{d+1} \,dxdt\Big)^{1/(d+1)}.
$$

By H\"older's inequality and in light of the fact that $1-\beta>0$
$$
\sum_{n=0}^{n_{0}}|u_{n}(0,x_{0})|
\leq N(d)(1-|x_{0}|)\Big(\sum_{n=0}^{n_{0}}
e^{n(1-\beta) (d+1)/d}\Big)^{d/(d+1)}
$$
$$
\times \Big( \int_{C}\psi^{\alpha}
f_{+}^{d+1}\sum_{n=0}^{n_{0}}\eta_{n}^{d+1}\,dxdt\Big)^{1/(d+1)}
$$
$$
\leq N(d,\alpha)(1-|x_{0}|)e^{n_{0}(1-\beta)}
\Big( \int_{C}\psi^{\alpha}
f_{+}^{d+1} \,dxdt\Big)^{1/(d+1)}
$$
\begin{equation}
                                                 \label{8.9.20}
\leq N(d,\alpha)\psi^{\beta}(x_{0})
\Big( \int_{C}\psi^{\alpha}
f_{+}^{d+1} \,dxdt\Big)^{1/(d+1)}.
\end{equation}

If   $n>n_{0}$, then by Lemma \ref{lemma 8.7.10}
$$
|u_{n}(0,x_{0})|\leq N(d)(1-s_{n}^{2})^{1/2}
\Big( \int_{C}\eta_{n}^{d+1}
f_{+}^{d+1}\,dxdt\Big)^{1/(d+1)}
$$
 $$
\leq N(d) e^{-n \beta}\big(\int_{C}  \psi^{\alpha}
\eta_{n}^{d+1}
f_{+}^{d+1}\,dxdt\Big)^{1/(d+1)}.
$$
Since $\beta>0$,
$$
\sum_{n>n_{0}}|u_{n}(0,x_{0})|\leq N(d)\Big(\sum_{n>n_{0}}
e^{-n \beta (d+1)/d)}\Big)^{d/(d+1) }
$$
$$
\times \Big(\int_{C}\psi^{\alpha}
f_{+}^{d+1}\,dxdt\Big)^{1/(d+1)}\leq N(d,\alpha)
e^{-n_{0} \beta }\Big(\int_{C}\psi^{\alpha}
\det D^{2}u \,dx\Big)^{1/(d+1)}
$$
$$
\leq N(d,\alpha)\psi^{\beta}(x_{0})
\Big(\int_{C}\psi^{\alpha}
\det D^{2}u \,dx\Big)^{1/(d+1)}.
$$
Upon combining this with \eqref{8.9.20} and Lemma
\ref{lemma 8.8.10} we get \eqref{8.5.10}.
The lemma is proved.

{\bf Proof of Theorem \ref{theorem 8.6.10}}. 
For small $\varepsilon>0$ set $C^{\varepsilon}
=[\varepsilon,T-\varepsilon)\times B_{1-\varepsilon}$
and
observe that $u\in W^{1,2}_{d+1}(C^{\varepsilon})$.
If the obvious version of \eqref{8.5.10} is true
with the objects with $\varepsilon$,
then setting $\varepsilon\downarrow0$ we obtain 
\eqref{8.5.10} as is. Hence, we may assume that
$u\in W^{1,2}_{d+1}(C)$. After that,
as usual, we may assume that $u$ is a smooth function.
Then take $\varepsilon>0$ and set $u_{\varepsilon}=
u-\varepsilon(\psi+T-t)$, so that $v_{\varepsilon}$
is strictly convex and strictly increasing.  Observe that on $\bar C$ 
$$
f_{\varepsilon}:=
\big(\partial_{t}u_{\varepsilon}\det D^{2}u_{\varepsilon}\big)^{1/(d+1)}
\geq\varepsilon 2^{d/(d+1)},
$$
and  $f_{\varepsilon}$ is smooth. We extend it
for $t\in(T,T+1]$, $x\in \bar B$, so that it remains nonnegative, smooth, becomes zero for $t=T+1$,
and
$$
\int_{(T,T+1)\times B}f^{d+1}_{\varepsilon}\,dxdt
\leq\varepsilon.
$$
Then define $C'=[0,T+1)\times B$
and  introduce $v_{\varepsilon}$ as a unique
$W^{1,2}_{\infty}(C')$-function which is convex
in $x$, increases in $t$, vanishes on $\partial'C'$
and satisfies
$$
\partial_{t}v_{\varepsilon}\det D^{2}
v_{\varepsilon}=f^{d+1}_{\varepsilon}
$$
in $C'$ (a.e.). In light of Lemma 
\ref{lemma 8.25.1}, to prove the theorem,
it suffices to prove that in $C$
\begin{equation}
                                      \label{8.29.1}
\bar u_{\varepsilon}:=u_{\varepsilon}-\inf_{\partial'C}\geq
v_{\varepsilon}
\end{equation}

Since $v_{\varepsilon}\leq0$, \eqref{8.29.1}
holds on $\partial'C$. Furthermore, both
$u_{\varepsilon}$ and $v_{\varepsilon}$
satisfy the same equation
$$
\inf_{(r,a)\in A}[r\partial_{t}w
+a^{ij}D_{ij}w-(d+1)f_{\varepsilon}\sqrt[d+1]{r\det a}]=0
$$
in $C$. Hence, \eqref{8.29.1} follows by the maximum
principle and the theorem is proved.

We now turn to estimates for parabolic operators of main type.
 \begin{theorem}  

                                           \label{theorem 8.11.10}

Let $ a(t,x)=(a^{ij}(t,x))$ and $r(t,x)$ be a $ d\times d$-symmetric
nonnegative definite matrix-valued measurable function and
a nonnegative measurable function on $C$, respectively,
such that $ r+\tr a>0$ in $ C$. Let $\alpha\in[0,(d+1)/2)$
 $u\in W^{1,2}_{d+1,\loc}(C)\cap C(\bar C)$.
Introduce
$$ 
L_{0}u=r\partial_{t}u+ a^{ij}D_{ij}u .
$$
Then, 
for any $ x_{0}\in B $, $(0/0:=0)$
$$
u(0,x_{0}) \leq \sup_{\partial'C} u
$$
\begin{equation}
                                         \label{8.9.70}
+N(d,\alpha)\psi^{\beta}(x_{0})
\Big(\int_{C}\psi^{\alpha}I_{Lu<0}
(r\det a)^{-1}|L_{0}u|^{d+1} \,dxdt\Big)^{1/(d+1)}.
\end{equation}
\end{theorem}

Proof. As in the proof of Theorem \ref{theorem 8.6.10},
we may assume that $u\in W^{1,2}_{d+1}(C)$.
Since $r+\tr a>0$, by the homogeneity of estimate
\eqref{8.9.70} we may also assume that
\begin{equation}
                                  \label{8.11.30}
r+\tr a\equiv 1.
\end{equation}

{\em A particular case\/}. Suppose that $r\det a\geq \varepsilon$,
where $\varepsilon>0$ is a constant.  Set
$$
f:=-\frac{1}{(d+1)}(r\det a)^{-1/(d+1)}L_{0}u
$$
and find a sequence of  $f_{n}\in L_{d+1}(C)$
such that $f_{n}\to f$ in $L_{d+1}(C)$,
$f_{n}$ are smooth, and vanish near
$\partial'C$. Then $(f_{n})_{+}^{d+1}
=\psi(\tilde f_{n})_{+}^{d+1}$, where
$\tilde f_{n}=f_{n}\psi^{-1/(d+1)}$
is a smooth function vanishing for $t=T$.

In that case, as we know, for each $n$,
there exists a unique
$W^{1,2}_{\infty}(C)$-function $v_{n}$ which is
concave  in $x$, decreases in $t$, vanishes  on $\partial' C$, and
satisfies
\begin{equation}
                                                \label{8.8.10}
-\partial_{t}v_{n}\det(- D^{2}v) =(f_{n})_{+}^{d+1}
\end{equation}
in $C$ (a.e.).  
We also know that $v_{n}$ also satisfies
\begin{equation}
                                                \label{8.11.50}
\sup_{(r,a)\in A}[r\partial_{t}v_{n}+a^{ij}D_{ij}v_{n}+f_{n}(d+1)\sqrt[d+1]{r\det a}]=0
\end{equation}
and for $u$ we obviously have
$$
\sup_{(r,a)\in A}[r\partial_{t}u+a^{ij}D_{ij}u+f(d+1)\sqrt[d+1]{r\det a}]
\geq0.
$$
It follows that for $w_{n}=u-\sup_{\partial'C}u-v_{n}$ we have
$$
\sup_{(r,a)\in A}[r\partial_{t}w_{n}+a^{ij}D_{ij}w_{n}+(f-f_{n})(d+1)\sqrt[d+1]{r\det a}]\geq0
$$
and there exists an $A$-valued function $(r(t,x),a(t,x))$
such that in $C$ (a.e.)
$$
r\partial_{t}w_{n}+a^{ij}D_{ij}w_{n}+(f -f_{n})(d+1)\sqrt[d+1]{r\det a}\geq0.
$$

By Theorem 3.2.2 of \cite{Kr_18}
$$
w_{n}\leq N\|f_{n}-f\|_{L_{d+1}(C)},
$$
where $N$ is independent of $n$. This and 
Theorem \ref{theorem 8.6.10}
show that
$$
u(0,x_{0})-\sup_{\partial'C}u\leq \nlimsup_{n\to\infty} v_{n}(0,x_{0})
\leq N(d,\alpha)\psi^{\beta}(x_{0})
\|\psi^{\alpha}f_{+}\|_{L_{d+1}(C)},
$$
which
completes considering this case.

{\em General case\/}. In order  to drop the additional assumptions
 we take $\gamma,\delta>0$,
 take $\psi$  from Lemma 3.1.8 of \cite{Kr_18}
  introduce  $ L^{\delta}= L+\delta 
( \partial_t+\Delta) $,
and apply the  above result  to  
$u^{\gamma}=u-\gamma(\psi+1)$ and  $ L ^{\delta}$.
Then following almost  word for word  
the appropriate parts of the proof
of Lemma 3.2.4 of \cite{Kr_18}, we  arrive at  \eqref{8.9.70}. 
The theorem is proved.

 \mysection{General parabolic equations in $C$}
                                      \label{section 8.11.20}

Here we generalize Theorem \ref{theorem 8.11.10}
for equation with lower order terms.
\begin{theorem}
                                         \label{theorem 8.11.20}
Under the assumptions of Theorem \ref{theorem 8.11.10}
let $b=(b^{i}(t,x))$ be $\bR^{d}$-valued measurable
function on $C$  and $c=c(t,x)$ be a nonnegative 
 measurable bounded
function on $C$ and let  $|b|\leq K\mu$ in $B$,
where $K$ is a fixed constant and $\mu=\mu(t,x)$
is the smallest eigenvalue of $a=a(t,x)$. Then
estimate \eqref{8.9.70} holds again with 
$N=N(d,\alpha,K)$  
if we replace $L_{0}u$ in \eqref{8.9.70} with
\begin{equation}
                                              \label{8.11.20}
Lu=r\partial_{t}u+a^{ij}D_{ij}u+b^{i}D_{i}u-cu,
\end{equation}
and $\sup_{\partial' C}u$ with $\sup_{\partial' C}u_{+}$.

\end{theorem}

 Proof. As in  the proof of 
Theorem \ref{theorem 8.11.10} we assume \eqref{8.11.30},
 $u\in W^{1,2}_{d+1}(C)$, 
and first also assume that  $r\det a\geq \varepsilon$,
where $\varepsilon>0$ is a constant. Set
$$
f:=-\frac{1}{d+1}(r\det a)^{-1/(d+1)}L u
$$
and define $f_{n}$ and $v_{n}$ in the same way
as in the proof Theorem \ref{theorem 8.11.10}.
 We also take the function $w$ from Lemma
\ref{lemma 8.11.3} with $\gamma=1-\beta$, take $N(d,\alpha)$ from
\eqref{8.5.10} and set  
$$
N(d,\alpha)
\psi^{\beta}(x_{0})\Big(\int_{C}\psi^{\alpha} 
(f_{n})_{+}^{d+1}\,dxdt\Big)^{1/(d+1)}=:N_{n}\psi^{ \beta}(x_{0}).
$$     

Observe that, as for any concave function vanishing on
$\partial B$, we have $|Dv_{n}|\leq v_{n}/(1-|x|)\leq 2v_{n}/\psi$.
Also note that for $\bar u=u-\sup_{\partial B}u_{+}$
we have
$
L\bar u\geq Lu$.
Then for the function
$$
\phi_{n}=v_{n}+2KN_{n}w
$$
owing to \eqref{8.5.10}, \eqref{8.11.4},
 and \eqref{8.11.50},   we get
$$
L\bar u\geq-(d+1)(r\det a)^{1/(d+1)}(f-f_{n})
+r\partial_{t}v_{n}+a^{ij}D_{ij}v_{n}
$$
$$
=:g_{n}+ Lv_{n}-b^{i}D_{i}v_{n}+cv_{n}\geq 
g_{n}+Lv_{n}-2|b| v_{n}/\psi
$$
$$
\geq g_{n}+Lv_{n}-2N_{n}|b| \psi^{-(1-\beta)}\geq
g_{n}+ L\phi_{n}.
$$
By Theorem 3.2.2 of \cite{Kr_18}
 $$
\bar u(0,x_{0})\leq\phi_{n}(0,x_{0})+N\|f-f_{n}\|_{L_{d+1}
(C)},
$$ 
where $N$  is independent of $n$.
After that it only remains to use 
Theorem \ref{theorem 8.24.1},
let $n\to\infty$, and recall that, as follows from
\eqref{9.17.1}, $w\leq N\psi^{\beta}$, where 
$N$ depends only on $d,\alpha,K$.
This proves the theorem in our
particular case.

The general case is dealt with as at the end of the proof
of Theorem \ref{theorem 8.11.10}.
The theorem is proved.

\mysection{The case of general smooth cylinders}
                                      \label{section 8.15.10}

Let $\Omega$ be a $C^{1,1}$ bounded domain in $\bR^{d}$,
$a=(a^{ij}(t,x))$ be $d\times d$ symmetric matrix-valued
function on $\bR^{d+1}$, $b=(b^{i}(t,x))$ be $\bR^{d}$-valued
function on $\bR^{d+1}$, and $c=c(t,x)$ be a nonnegative bounded
function on $\bR^{d+1}$.
Suppose that these functions are measurable.
 We fix two numbers $\delta,K >0$
and assume that, for all values of arguments and $\lambda\in\bR^{d}$,
\begin{equation}
                                                        \label{8.15.50}
\delta^{-1}|\lambda|^{2}\geq
a^{ij}\lambda^{i}\lambda^{j}\geq \delta|\lambda|^{2},
\quad |b|\leq K, \quad  c\geq 0.
\end{equation}
Introduce $\Pi=[0,T)\times\Omega$, $\partial'\Pi=\partial \Pi
\setminus (\{0\}\times \bar\Omega)$,
$$
Lu=\partial_{t}u+a^{ij}D_{ij}u+b^{i}D_{i}-cu.
$$
Finally, take $\alpha\in[0,(d+1)/2)$ and define
$\rho_{\Omega}(x)=\dist(x,\Omega^{c})$. 

\begin{theorem}
                                 \label{theorem 8.15.10}

Let $u\in  W^{1,2}_{d+1,\loc}(\Pi)\cap C(\bar\Pi)$.
Then, 
for any $(t, x) \in \Pi $, 
\begin{equation}
                            \label{8.15.10}
u(t,x ) \leq \sup_{\partial'\Pi}u_{+}+N \rho_{\Omega}^{\beta}(x )
\|\rho_{\Omega}^{\alpha} 
  (Lu)_{-}\|_{L_{d+1}(\Pi)},
\end{equation}
where 
the constant $N$ depends only on $\Omega$, $\delta$, and $K$.
\end{theorem}

To prove this theorem we need three lemmas.
 
\begin{lemma}
                 \label{lemma 8.11.10}
Set $B^{+}_{r}=B_{r}\cap\{x\in\bR^{d}:x^{1}>0\}$, $B^{+}=B^{+}_{1}$,
$C^{+}=[0,T)\times B^{+}$,
$\partial'C^{+}=\partial C^{+}\setminus(
\{0\}\times\bar B^{+})$
and suppose that $u\in W^{1,2}_{d+1,\loc}(C^{+})\cap C(\bar
C^{+})$,
$u\leq 0$ on $\partial' C^{+}$. Also assume that
 the coefficients of $L$ are infinitely differentiable
and $f:=Lu \in L_{d+1}(C^{+})$.
 Then in $C^{+}$    
we have
\begin{equation}
                                   \label{8.15.40}
u(t,x ) \leq N (x^{1})^{\beta}
\|M^{\alpha} 
  f_{-}\|_{L_{d+1}(C^{+})},
\end{equation}
where $M^{\alpha}(x)=(x^{1})^{\alpha}$  and the constant
$N$ depends only on $d,\delta$, and $K$.

\end{lemma}

Proof. Find
a $C^{1,1}$ domain $B'$ which contains $B^{+}_{2}$,
  has $B_{2}\cap\{x\in\bR^{d}:x^{1}=0\}$
as part of its boundary, and is $C^{1,1}$-diffeomorphic
to $B$. Set $C'=(0,T)\times B'$.

Then, extend $f$ as zero outside $C^{+}$
and denote by $v$ the 
function of class $W ^{1,2}_{d+1}(C')$ satisfying
\begin{equation}
                                \label{8.15.60}
 Lv=-f_{-}
\end{equation}
in $C'$ (a.e.) and vanishing on
$\partial' C'$. It is a classical fact that
such a function exists and is unique.
By the maximum principle $v\geq 0$, $v\geq u$ in $C^{+} $.
Furthermore, is we apply the diffeomorphism
mentioned above to \eqref{8.15.60} then we will see
that the image $v'$ of $v$ will satisfy the equation
$$
 L'v'=-f'_{-}
$$
in $C$ (a.e.),
where $f'_{-}$ is the image of $f_{-}$ and $ L'$
is the image of $L$.

By Theorem \ref{theorem 8.11.20}
\begin{equation}
                                                        \label{8.15.70}
v'(x')\leq N\psi^{\beta}(x')
\|\psi^{\alpha/d} f'_{-}\|_{L_{d+1}(C )},
\end{equation}
where $\psi(x')$ is the distance from $x'$ to the 
boundary of $B$. If $x'$ is the image of $x$,
$\psi(x')$ is   comparable to the distance of $x$ to
the boundary of $B'$, since the diffeomorphism
and its inverse are   Lipschitz continuous.

It only remains to write down
\eqref{8.15.70} in the original coordinates and
use the fact that $v\geq u$  and that for $x\in B^{+}$
its distance  to
the boundary of $B'$ equals $x^{1}$.
The lemma is proved.

Next we use a well-known fact
(see, for instance, Lemma 8.8
in \cite{Kr_11})
that  there exists a function $\Psi\in C^{1,1}(\bar \Omega)$
such that, for a constant $N$ depending only on
$\Omega$, $\delta$, and $K$, we have in $\Omega$
(for any $t$)
\begin{equation}
                               \label{8.15.20}
 N\rho_{\Omega}\geq \Psi , \quad \Psi_ \geq \rho_{\Omega},
\end{equation}
\begin{equation}
                                \label{8.15.30}
L\Psi +c\Psi \leq-1.
\end{equation}

\begin{lemma}
                         \label{lemma 8.15.10}
Let the coefficients of $L$ be infinitely differentiable.
Take $\rho_{0}\in(0,1)$, $x_{0}\in \Omega$ and suppose that
$2\rho:=\rho_{\Omega}(x_{0})\geq  \rho_{0}$. Let $\gamma\in (0,1]$ and let
$\Phi$ be the classical bounded
 solution of $L\Phi=0$ in $[0,\infty)\times B_{\rho }(x_{0})$
with boundary condition $\Phi=\Psi^{\gamma}$ for 
$|x-x_{0}|=\rho$. Then
$$
\Phi(0,x_{0})\leq[1-\varepsilon \rho_{0}]^{\gamma} \Psi^{\gamma}(x_{0}),
$$
where $\varepsilon>0$ depends only on $K$, $\delta$, and $\Omega$.
\end{lemma}
 
Proof. A simple argument based on the maximum principle shows
that it suffices to concentrate on $\gamma=1$. 
In that case first we note that 
by the maximum principle $\Phi\geq0$. Therefore, $(L+c)\Phi\geq0$
and
$v:=\Psi-\Phi$ satisfies $(L+c)v\leq-1$ and $v=0$
for 
$|x-x_{0}|=\rho$. Then elementary barriers 
show that $v(t,x_{0})\geq \varepsilon \rho ^{2}$,
where $\varepsilon>0$ depends only on $K$, $\delta$,
and the diameter of $\Omega$.
Since $\rho_{\Omega}(x_{0})$ and $\Psi(x_{0})$
are comparable we have (we use $\varepsilon$ as a generic
constant $>0$ depending only on $K$, $\delta$, and $\Omega$)
$$
\rho^{2}=\rho(1/2)\rho_{\Omega}(x_{0})\geq
\varepsilon \rho_{0}\Psi(x_{0}),
$$
$$
v(0,x_{0})=
\Psi(x_{0})-\Phi(0,x_{0})\geq \varepsilon \rho_{0}\Psi(x_{0}),
$$
and the lemma is proved.

Below by $\rho_{0}$ we mean a number $>0$
such that  any point $\bar x_{0}\in\partial\Omega$
is the only common point of $\partial\Omega$
and the closure of a ball, say $B_{\rho_{0}}(y_{0})$,
 belonging to $\Omega$
with radius $\rho_{0}$.
Since $\Omega\in C^{1,1}$ such $\rho_{0}>0$ exists.
We further decrease $\rho_{0}$, if necessary, so that
there is a $C^{1,1}$-diffeomorphism with its first- and second-order derivatives
and the first- and second-order derivatives of its inverse
bounded by a constant depending only on $\Omega$ and
mapping
$B_{2\rho_{0}}(\bar x_{0})\cap \Omega$ onto
$B^{+}_{1}$ and $\bar B_{2\rho_{0}}(\bar x_{0})\cap 
\partial\Omega$ onto $\bar B^{+}_{1}\cap\{x^{1}=0\}$.

In the following lemma we consider the case that
$ \rho_{\Omega}(x_{0})\leq  \rho_{0}$ and denote
by $\bar x_{0}$ a   point on $\partial \Omega$
such that $\rho_{\Omega}(x_{0})=|\bar x_{0}-x_{0}|$.
By assumption there exists $y_{0}\in\Omega$ such that
$$
\bar B_{\rho_{0}}(y_{0})\cap\partial\Omega=\{\bar x_{0}\},
$$
and since both balls $B_{\rho_{0}}(y_{0})$ and
$B_{\rho_{\Omega}(x_{0})}(x_{0})$ touch $\partial\Omega$
at $\bar x_{0}$ and the former ball has smaller radius,
the points $\bar x_{0},x_{0}$, and $y_{0}$ lie
on the same line and
\begin{equation}
                                                \label{8.16.10}
\dist(x_{0},B^{c}_{\rho_{0}}(y_{0}))=\rho_{\Omega}(x_{0}).
\end{equation}

\begin{lemma}
                                            \label{lemma 8.15.20}
Let the coefficients of $L$ be infinitely differentiable.
Take $\rho_{0}\in(0,1)$, $x_{0}\in \Omega$ and suppose that
$\rho:= \rho_{\Omega}(x_{0})\leq  \rho_{0}$. 
Take $y_{0}$ introduced before the lemma.
Let $\gamma\in (0,1]$ and let
$\Phi$ be the classical bounded
 solution of $L\Phi=0$ in $(0,\infty)\times B_{ \rho_{0} }(y_{0})$
with boundary condition $\Phi=\Psi^{\gamma}$ on 
$(0,\infty)\times \partial B_{ \rho_{0} }(y_{0})$. Then
$$
\Phi(0,x_{0})\leq[1-\varepsilon \rho_{0}]^{\gamma} \Psi^{\gamma}(x_{0}),
$$
where $\varepsilon>0$ depends only on $K$, $\delta$, and $\Omega$.
\end{lemma}

Proof. As in the proof of Lemma \ref{lemma 8.15.10},
we concentrate on the  case that $\gamma=1$ and
we have $(L+c)v\leq-1$. Simple barriers show that
(recall that $\rho_{0}$ is fixed)
$v(t,x)\geq \varepsilon(\rho_{0}-|x-y_{0}|)$ in $B_{ \rho_{0} }(y_{0})$,
where $\varepsilon>0$ depends only on $d$, $\delta$, $K$, and
$\rho_{0}$. Owing to \eqref{8.16.10} this shows that
$v(0,x_{0})\geq \varepsilon \rho_{\Omega}(x_{0})$
and we are done because $\rho_{\Omega}(x_{0})$ and $\Psi(x_{0})$
are comparable. The lemma is proved.

{\bf Proof of Theorem \ref{theorem 8.15.10}}.
We can certainly concentrate on the case that
$u$ and the coefficients of $L$ are infinitely differentiable
in $\bar \Pi$. In that case the function
$$
v=\frac{u}{\Psi^{\beta}}
$$
is continuous in $\bar \Pi$, equals zero on $\partial'\Pi$ and, hence,
attains its maximum value  in $\bar \Pi$ at a point $(t_{0},x_{0})\in \Pi$:
\begin{equation}
                                                      \label{8.16.30}
M:=\frac{u(t_{0},x_{0})}{\Psi^{\beta}(x_{0})}\geq \frac{u(t,x)}{\Psi^{\beta}(x)}
\quad\forall (t,x)\in \Pi.
\end{equation}
If $M$ is less than zero, we have nothing
to prove. Therefore, we assume that
$$
u(t_{0},x_{0})>0
$$ 
and consider two cases:

(a) $\rho_{\Omega}(x_{0})\geq \rho_{0}$,

(b) $\rho_{\Omega}(x_{0})< \rho_{0}$.

In case (a), in $C':=[t_{0},T)\times B_{\rho_{\Omega}(x_{0})/2}(x_{0})$ we have $u=v+h$,
where $v$ is the classical solution of $Lv=Lu$ 
in $C'$ with zero boundary value and $h$ is the classical
solution of $Lh=0$ in $C'$  with boundary condition
$h=u$ on $\partial' C'$. We apply the Aleksandrov
estimate to $v$ and take into account that
on $C'$, $\rho_{\Omega}$ and $\Psi$ are comparable to constant one.
Then we see that
\begin{equation}
                                                      \label{8.16.20}
v(t_{0},x_{0})\leq N\| 
  (Lu)_{-}\|_{L_{d+1}(C')}\leq
N \Psi^{\beta}(x_{0} )
\|\rho_{\Omega}^{\alpha} 
  (Lu)_{-}\|_{L_{d+1}(\Pi)}.
\end{equation}

In what concerns $h$, observe that owing to \eqref{8.16.30}
and the fact that $u(T,x)=0$ on $\Omega$, 
by the maximum principle,
 $h$ is less
than the bounded classical solution $\Phi$ of the equation $L\Phi=0$ in $[t_{0},\infty)\times B_{\rho_{\Omega}(x_{0})/2}(x_{0})$
with boundary condition $M\Psi^{\beta}$. By Lemma \ref{lemma 8.15.10},
$h(t_{0},x_{0})\leq \Phi(t_{0},x_{0})\leq M\varepsilon\Psi^{\beta}(x_{0})$,
where $\varepsilon\in(0,1)$ depends only on $K,\delta,\alpha$,
and $\Omega$. It follows that
$$
M\leq N\|\rho_{\Omega}^{\alpha} 
  (Lu)_{-}\|_{L_{d+1}(\Pi)}+M\varepsilon,\quad
M\leq N\|\rho_{\Omega}^{\alpha} 
  (Lu)_{-}\|_{L_{d+1}(\Pi)},
$$
which yields \eqref{8.15.10}.

In case (b), take the ball $B_{\rho_{0}}(y_{0})$
introduced before Lemma \ref{lemma 8.15.20} 
and set $C'=[t_{0},T)\times B_{\rho_{0}}(y_{0})$.
Then in $C'$
we have
$u=v+h$, where $v$ is the solution of $Lv=Lu$
in $C'$ with zero boundary condition
and $h$ satisfies $Lh=0$ and equals $u$ on
$\partial' C'$. As in case (a)
by the maximum principle and Lemma \ref{lemma 8.15.20}
we have $h(t_{0},x_{0})\leq M\varepsilon \Psi^{\beta}(x_{0})$.

In what concerns $v$ observe that in $C' $
 by the maximum
principle it is less than $w$ defined as
$W^{2}_{d+1,\loc}(C'')\cap C(\bar C'')$-solution of
$Lw=-f_{-}$ in
$$
  C''=[t_{0},T)\times(B_{2\rho_{0}}(\bar x_{0})\cap \Omega)
$$   
vanishing on its parabolic  boundary, where
$$
f:=I_{C'}Lu.
$$
By an argument similar to the one used in the proof of Lemma
\ref{lemma 8.11.10} we obtain that
$$
w(t_{0},x_{0})\leq N\rho^{\beta}_{\Omega}(x_{0})\|\rho_{\Omega}^{\alpha}f_{-}
\|_{L_{d+1}(C'')}
\leq N\rho^{\beta}_{\Omega}(x_{0})\|\rho_{\Omega}^{\alpha}f_{-}
\|_{L_{d+1}(\Pi)}.
$$
Since $\rho_{\Omega}$ and $\Psi$ are comparable, we conclude
$$
M\leq M\varepsilon+N\|\rho_{\Omega}^{\alpha}f_{-}
\|_{L_{d+1}(\Pi)}
$$
and this proves the theorem.

\mysection{Estimates for stochastic integrals}
                                                   \label{section 9.1.10}

Let $(\Omega,\cF,P)$ be a complete probability
space, $\{\cF_{t},t\geq0\}$ be an increasing filtration
of $\sigma$-fields $\cF_{t}\subset\cF$ each of 
which is complete relative to $\cF,P$. Suppose
that on $(\Omega,\cF,P)$ we are given a $d_{1}$-dimensional Wiener process which is $\cF_{t}$-adapted
and such that $w_{t}-w_{s}$ are independent of $\cF_{s}$
as long as $0\leq s\leq t<\infty$.

Fix some constants $\delta\in(0,1]$ and $K\geq0$.
Let $\sigma_{t}=\sigma_{t}(\omega)$ be a progressively
measurable with respect to $\{\cF_{t},t\geq0\}$,
$d\times d_{1}$-matrix valued process such that
$$
\delta^{-1}|\lambda|^{2}\geq a_{t}^{ij}\lambda^{i}\lambda^{j}
\geq\delta|\lambda|^{2}
$$
for all $\lambda\in\bR^{d}$, $t\geq0$, and
$\omega\in\Omega$, where $a=(1/2)\sigma\sigma^{*}$.

Let $b_{t} $ be a progressively
measurable with respect to $\{\cF_{t},t\geq0\}$,
$\bR^{d}$-valued process such that
$|b_{t}|\leq K$
for all  $t\geq0$ and
$\omega\in\Omega$.

Introduce
$$
x_{t}=\int_{0}^{t}\sigma_{s}\,dw_{s}+
\int_{0}^{s}b_{s}\,ds.
$$

Let $G$ be a bounded $C^{1,1}$ domain in $\bR^{d}$
containing the origin, set
$$
\rho_{G}(x)=\dist(x,G^{c}),
$$
take a nonnegative Borel function $f$ on $
(0,T)\times G$,
a number
$$
\alpha\in[0,(d+1)/2),
$$
and set
$$
\beta=\frac{1}{2}-\frac{\alpha}{d+1}.
$$
Here is the result of this section.

\begin{theorem}
                                                  \label{theorem 9.1.10}
There exists a constant $N$, depending only on $G$, $\delta$, $\alpha$, and $K$, such that
\begin{equation}
                     \label{9.1.20}
E\int_{0}^{\tau}f(t,x_{t})\,dt
\leq N\rho_{G}^{\beta}(0)
\Big(\int_{(0,\infty)\times G}
\rho_{G}^{\alpha}(x)f^{d+1}(t,x)\,dxdt\Big)^{1/(d+1)},
\end{equation}
where $\tau$ is the first exit time of $x_{t}$ from $G$.
\end{theorem}

Proof. A usual measure-theoretic argument
shows that it suffices to prove \eqref{9.1.20}
for bounded $f$ vanishing for $t\geq T$
with arbitrary $T\in(0,\infty)$.
Let $A$ be the collection of couples
$(a,b)$, where $a$ is a symmetric $d\times d$
matrix and $b\in \bR^{d}$ such that
$$
|b|\leq K,\quad \delta^{-1}|\lambda|^{2}\geq a ^{ij}\lambda^{i}\lambda^{j}
\geq\delta|\lambda|^{2}
$$
for all $\lambda\in\bR^{d}$.

As it follows from   \cite{Kr_18}
there exists a unique $u\in W^{1,2}_{d+1}
((0,T)\times G)$
vanishing on the parabolic boundary
of $(0,T)\times G$ and satisfying
$$
\sup_{(a,b)\in A}[\partial_{t}u+a^{ij}D_{ij}u+b^{i}D_{i}u+f]=0
$$
in $(0,T)\times G$ (a.e.). By It\^o's formula
(see, for instance, \cite{Kr_77})
$$
u(0)=-E\int_{0}^{\tau\wedge T}[\partial_{t}u(t,x_{t})+a^{ij}_{t}D_{ij}u(t,x_{t})+b^{i}_{t}
D_{i}u(t,x_{t})]\,dt
$$
$$
\geq
E\int_{0}^{\tau\wedge T}f(t,x_{t})\,dt
=E\int_{0}^{\tau}f(t,x_{t})\,dt.
$$

After that it only remains to 
apply Theorem \ref{theorem 8.15.10} first
observing that
there exist a measurable $A$-valued function
 $(a(x),b(x))$ such that
$$
\partial_{t}u(t,x)+a^{ij}(t,x)D_{ij}u(t,x)+b^{i}(t,x)D_{i}u(t,x)+f(t,x)=0
$$
in $(0,T)\times G$ (a.e.).  The theorem is proved.

If Borel $\Gamma\subset (0,\infty)\times G$, then
$$
G(\Gamma)=E\int_{0}^{\tau}I_{\Gamma}(t,x_{t})\,dt
$$
(the so-called Green's measure) is the mean time that the trajectory
$(t,x_{t}),t\in[0,\tau]$, spends in $\Gamma$,
the mean time it occupies $\Gamma$
before time $\tau$.
By taking $f=I_{\Gamma}$ in \eqref{9.1.20}
we come to the following.

\begin{corollary}
We have
$$
G(\Gamma)\leq 
N\rho_{G}^{\beta}(0)
\Big(\int_{\Gamma}\rho_{G}^{\alpha} (x)\,dxdt\Big)^{1/d}.
$$
\end{corollary}

{\bf Acknowledgement}. The article was finished
during the author's visit to the University
of Bielefeld on the invitation of Michael Roeckner.
This is greatly appreciated. The author is also sincerely grateful
to A.I. Nasarov for providing important information
on the subject of the paper.

\end{document}